\documentclass{article}

\newtheorem{thm}{Theorem}[section] %

\newcommand{\eqn}{\begin{eqnarray}}

\newcommand{\eeqn}{\end{eqnarray}}

\begin{document}

\title{Eigenvalues of Transmission Graph Laplacians}

\author{Sylvain E. Cappell and Edward Y. Miller  }

\date{}

\maketitle

\begin{flushleft}

\footnotesize

\centerline{ABSTRACT}
The standard notion of the Laplacian of a graph is generalized to the
setting of a graph with the extra structure of a ``transmission'' system.
A transmission system is a mathematical representation of a means of transmitting (multi-parameter) data along directed
edges from vertex to vertex. The associated transmission graph  Laplacian
is shown to have many of the former properties of the classical
case, including:  an upper Cheeger type bound on the second eigenvalue
minus the first of a geometric isoperimetric
character, relations of this difference of eigenvalues to diameters
for $k$-regular graphs, eigenvalues for Cayley graphs with transmission
systems. An especially natural transmission system arises in the context of
a graph endowed with an association. Other relations to  transmission systems arising naturally in
quantum mechanics, where the transmission matrices are  scattering matrices, are made.
As a natural merging of graph theory and matrix theory, there are numerous
potential applications,  for example to random graphs and random matrices.

\vspace{.1in}

\S 1: Introduction.

\S 2: Definition of the Transmission Graph Laplacian and the Generalized Cheeger Estimate.

\S 3: A Symmetric Graph Transmission System for Associations.

\S 4: Quantum Mechanical Example of Transmission Systems: Invertible and (Hermitian) Symmetric.

\S 5: The Dual of a Graph and its Transmission Systems.

\S 6: Estimates 1: Range of Values and Cheeger's Upper Bound Generalized.

\S 7: Estimates 2: Eigenvalues and Diameters.

\S 8: Reinterpretation in terms of Capacities. Cheeger Constants.

\S 9: Cayley graphs. Their Eigenvalues.

\S 10: Eigenvalues for Graph Collapses and Amalgamations.

\S 11: A Categorical Approach. The push forward collapse.

\S 12: Morse Theory, Riemann Surfaces and Transmission systems.

\end{flushleft}

\section{Introduction.}  \label{section1}

In this paper the standard notion of the  Laplacian of a graph, or graph Laplacian,
  is generalized to the setting of
a graph $X$ with the added structure of a ``transmission system''.
A model to keep in mind is that at each vertex, say $v$, data is recorded
by a list of (complex) numbers. The length of the list is called the ``bit-rank''
since it is exactly the number of bits if the field $C$ were replaced by the
field $Z_2= \{0,1\}$ of two elements. This list is in the form of  a column vector, say  $a$,  in a vector
space of some fixed $C^{n(v)}$ to be regarded as the data bank at $v$ of bit-rank $n(v)$.
A ``transmission system''
is just a specification for each directed edge $\vec{e}$ of the graph of an associated
  matrix $P(\vec{e})$.
The matrix $P(\vec{e})$
for a directed edge $\vec{e}$ with  initial point $v$ and  terminal point $w$
may be  regarded as
 a  means of  transmitting
 data kept at the initial node $v$ recorded by a column vector, say $a \in C^{n(v)}$,
to data kept  at the terminal node $w$ in $C^{n(w)}$  by sending $ a \mapsto P(\vec{e}) \ a $.
Thus the sole condition on a general transmission matrix $P(\vec{e})$ is that it
be a $ n(w) \times n(v)$ matrix otherwise completely arbitrary.
 In the case that all the ranks are one  for all vertices  and all  mappings  the identity,
  the old notions of graph theory are recaptured.

 In section \ref{sectstart} the definition of the (normalized) graph Laplacian
 as defined by Chung \cite{Chung1} is reviewed along with the statement of Cheeger's
 geometric estimate for the second eigenvalue of the graph Laplacian. Patterned
 on this discussion, the transmission graph Laplacian for a transmission system
 is defined, its range of eigenvalues is clarified, Theorem \ref{thmrange}, and
  generalizations of the differential  geometric estimate of Cheeeger stated in
 Theorems \ref{thmCheegerupper1},\ref{thmCheegerupper2}.

 If $op(\vec{e})$ is the directed edge $\vec{e}$
 with the direction reversed to its opposite,  one says that the transmission system $P$
 is
 $$
 \begin{array}{ll}
(Hermitian) \ symmetric \ if &    P(op(\vec{e})) = P(\vec{e})^\star \ \\
\hspace{-0.4in} and \ is \\
invertible \ if  &  P(op(\vec{e})) = P(\vec{e})^{-1} \ . \\
\end{array}
$$
In the first case, the associated transmission graph Laplacian is
self-adjoint, so has real eigenvalues. In the second, the transmission
matrices are square and the transmission
system is equivalent to the topological notion of a local coefficient system, or flat vector bundle,
and the transmission graph Laplacian is just the canonical coupling of
this local system to the standard graph Laplacian\footnote{In the differential geometric setting, Laplacians of
 couplings  to a flat system which is  not necessarily unitary
were studied by the present authors in \cite{CappellMiller}.}. In the case that
the transmission system is both (Hermitian) symmetric and invertible,
each of $P(\vec{e})$ is unitary, so one says the transmission system
is strictly unitary in that case.

In sections \ref{sectassociation},\ref{sectquantum}, we exhibit and study some varied examples of symmetric transmission
systems arising in social and physical networks.

When graphs are used to represent (e.g., social) networks it is natural to consider assignments to
each vertex of a set of associated vertices, which for mathematical purposes could be arbitrary.
[For example, if the network is the social  network of friends, adjacency equaling  friendship, one might for
each person take the set of, say  face book friends of that person, or the set of backers of the same
sport team as that person, etc.] Such assignments will be called an association.
In studying such settings it's further natural to measure the the cohesion between two disjoint subsets, say $A,B$,
of vertices in terms of the total number of common associates for the pairs of vertices, $u,v$, $u$ in $A$ and
$v$ in $B$, with $u,v$ adjacent. Section \ref{sectassociation} studies such problems by means of the
eigenvalues of the Laplacian of  naturally arising symmetric transmission system defined by this association.
For the convenience of the reader this section is made free standing and the proofs, modeled
on the methods  of Chung \cite{Chung1}, serve as a pattern for the general graph transmission cases.

Another natural setting for associations is that in which there is second  graph on the same vertices,
perhaps recording other data or  relations. In this context, one may
 associate to a point  of the first graph  the  set of points which are at most a specified distance from that point
with respect to the adjacency relation of the second graph.
Note also that one can easily introduce an added density function on the set associated to a vertex
and incorporate this into the definition of the Laplacian also. This refinement is not dealt with here.

In section \ref{sectquantum} physically  arising  examples of such transmission systems
 are given in the setting of   the quantum mechanical transmission and reflection
 of beams by  potential barriers erected along the edges of a graph.
The scattering matrices for propagation along directed edge defines an
invertible transmission system, but not (Hermitian) symmetric.
 Surprisingly, based on quantum mechanical principle of  conservation of
probability, multiplying these by a fixed matrix yields a (Hermitian) symmetric
transition system, and so a self-adjoint operator which is here called
the quantum transmission graph Laplacian. (It differs in character from the standard
quantum graph Laplacian in the literature \cite{Gnutzman}.)

Another class of examples arises for the dual of a graph, see section \ref{sectdual},
which records the transmission of data through the vertices.

In this paper, it  is shown that
some of the  most important  properties of the eigenvalues of the graph Laplacian
 extend naturally to the transmission system  context and its associated transmission
 graph Laplacian. These are:
\begin{itemize}
\item The range of values of the eigenvalues of the
 graph Laplacian.
\item The geometric estimate of Cheeger   bounding the second eigenvalue
from above  in terms of splitting
$X$ into two pieces.
 This geometric estimate is  an isoperimetric
inequality, similarly for its generalization.
 \item the result of Nilli \cite{Nilli1,Chung2}
bounding the second  eigenvalue from below in terms of the diameter of the graph.
This estimate suggest the appropriate notion of ``Ramanujan graphs''.

 \item The relation of eigenvalues under a collapsing operation
on graphs \cite{Chung1}
\item A formula for the eigenvalues of Cayley graphs \cite{Murty}.
\end{itemize}

A more categorical approach to transmission systems appears in section
\ref{sectPushForward}. There are natural settings where the vector
space associated to a vertex does not come  equipped with a basis,
so more elegantly one assigns to each vertex a vector space
and appropriate linear mappings between them for each directed edge.
An application occurs in that section and also the next which is on
graphs embedded in Riemann surfaces.

It is natural to make many conjectures about the properties of the
transmission graph Laplacian. For example, the various improvements
of the result of Nilli might have analogues. Clearly, it is a  challenge
to attempt to  generalize all known results on the standard graph Laplacian
to this setting. In another direction, one may consider random choices of
strictly unitary
transmission systems and attempt to generalize theorems about
the eigenvalues of random unitary matrices to this context. Also,
one may randomize over $k$-regular graphs and strictly unitary transmission systems.
One may regard modifying the  transition matrices as an introduction of noise
into a system, or as the result of a spying operation at the edges and
study the implications. All these and more deserve additional study.

A primary foundational item not addressed  here is the question of the existence of
a  ``Ramanujan graph'' suitably defined in this extended setting,see the work of  Sarnak, Lubotzsky,
the survey of Murty, and the forthcoming book of Janwa and Rangachari \cite{Sarnak,Lubotzky,Murty, Janwa-Rangachari},
In this direction, the results of Friedman \cite{Friedman} suggest that random $k$-regular graphs
with random strictly unitary transmission systems may well satisfy the appropriate estimates.
Such an example would generalize notions of
$c$-expanders and  expander codes \cite{Sipser}.

\vspace{.3in}

In the classical case, the lowest eigenvalue of the graph Laplacian is zero and the associated
eigenvector is
$$
f(v) = \sqrt{d_v} \ for \ v \ a \ vertex \ of \ X
$$
where $d_v$ is the degree of the vertex, see section \ref{sectstart}.
In that setting, for a subset   $A$ of the vertices, $V(X)$ of the graph $X$, the natural volume is
the sum of the degrees of vertices in $A$.
$$
vol(A) =  \Sigma_{v  \in  A}  \ \ d_v =  \Sigma_{v \in  A} \ \  |f_1(v)|^2
$$

For $P$ a (Hermitian) symmetric transmission system, the transmission graph Laplacian
is self-adjoint and its eigenvalues real, so one may chose an eigenvector, say $f_1(v)$
of the associated transmission graph Laplacian with the least real eigenvalue, say $\lambda_1(X,P)$.
It will be seen that the correct generalization of the above classical volume of a set of vertices is
$$
vol_P(A) = \Sigma_{v  \in  V(X)} \ \ |f_1(v)|^2 \ .
$$
It is in these terms that the Cheeger geometric upper estimates is generalized, Theorems \ref{thmCheegerupper1}, \ref{thmCheegerupper2}.
The estimate is a geometrical upper estimate for the difference
$$
\lambda_2(X,P) - \lambda_1(X,P)
$$
of the second smallest eigenvalue minus the first.
This  geometric flavor is especially compelling in view of potential applications.

It is natural to think that the randomizing over strictly unitary transmission systems
will produce interesting distributions for $min_{ v \ in \ V(X) \  with \ f(v) \neq 0} \ \frac{|f(v)|}{\sqrt{d_v}}$
and  $max_{ v \ in \ V(X)} \ \frac{|f(v)|}{\sqrt{d_v}}$. These quantities also appear geometrically in  Theorem
\ref{thmCheegerupper2}.

This paper uses the exposition of Murty \cite{Murty} which provides
an elegant and direct approach to the properties of the eigenvalues
of the graph Laplacian in the case of $k$-regular graphs. In that case,
this study is equivalent to the study of the eigenvalues of the
adjacency matrix, as explained in section \ref{sectstart}.
This paper also uses   the beautiful book of Chung \cite{Chung1}.

It is a pleasure to acknowledge discussions with  Professor
Heeralal Janwa towards further understanding these new graph phenomena
and their potential applications.

\vspace{.3in}
\section{Definition of the Transmission Graph Laplacian and the Generalized Cheeger Estimate.} \label{sectstart}

We begin with a brief review  of the classical definition of the Graph Laplacian and some of its properties:
In this paper, the graph Laplacian of a graph $X$ is the (normalized) graph Laplacian as
defined by Chung in her book \cite{Chung1}. More explicitly, let the vertices of the graph $X$
be denoted by $V(X)$ and its edges by $E(X)$. As here an  emphasis will be on transmitting data
along edges, let $E^+(X)$ denote the set of directed edges of $X$, two for each edge. For
a directed edge, say $\vec{e} \in E^+(X)$, let the initial point of $\vec{e}$ be denoted by
$i(\vec{e})$ and the terminal point by $t(\vec{e})$. Let $op(\vec{e})$ denote the same edge
with the reversed, or opposite direction from that of $\vec{e}$. For a vertex $w$ of $X$, define the
degree of $w$, $deg(w) = d_w$ as  the number of  directed edges $\vec{e}$ with terminal point $w$,
$$
deg(w) = d_w = \Sigma_{ \vec{e} \ with \  t(\vec{e}) = w} \ \ (+1)
$$
This counts the (undirected) loops based at $w$ with multiplicity 2,  For a finite set $S$, let $|S|$
denote the cardinality of $S$. Let $n = |V(X)| $ be the number of vertices of the graph $X$.

Now to be as explicit as possible, order the vertices of $X$, say \newline
$V(X) = \{ v_1,v_2,\cdots, v_n\}$, and define the
adjacency matrix $A$ as the $n$ by $n$ matrix with $(i,j)$ entry defined
by the sum over directed edges with initial point $v_i$ and terminal point $v_j$,
$$
A_{i,j} = | \ \{ \vec{e} \ | \  i(\vec{e}) = v_i,\ t(\vec{e})= v_j \} \ |
$$
Define the normalized adjacency matrix via
$$
N(A)_{i,j} =  \frac{1}{\sqrt{ d_{v_i} \ d_{v_j} }}  \   |\  \{ \vec{e} \ | \  i(\vec{e}) = v_i,\ t(\vec{e})= v_j \} \ |
$$
Let $T$ be diagonal matrix which has  $(i,i)$  entry equal to the degree $d_{v_i}$ of the vertex $v_i$, and denote by
$T^{m}$ the diagonal matrix with  $(i,i)$ entry equal to  $(d_{v_i})^m$.

In these terms, the (normalized) graph Laplacian is the $n$ by $n$ matrix
$$
\Delta = Id - N(A) = T^{-1/2} \ ( T - A) T^{-1/2}
$$

In the case that the graph $X$ is $k$-regular, that is $d_v = k$ for all vertices $v$,
$\Delta = Id - (1/k) A = (1/k) (k \ Id - A)$ so the eigenvalues of $\Delta$ are exactly $1- (1/k) \hat{\lambda}(A)$
where $\hat{\lambda}(A)$ ranges over the eigenvalues of the adjacency matrix $A$.
[In the restricted setting of $k$-regular graphs it is the matrix $k Id - A$ which is sometimes
referred to as the graph Laplacian; this is not done here.]
\vspace{.3in}

The salient features of the (normalized) graph Laplacian $\Delta$ are \cite{Chung1}:
\begin{itemize}
\item $\Delta$ is real symmetric, so has real eigenvalues. These lie in the range $0$ to $+2$.
Listed by increasing magnitude let them be
$$
0 \le \lambda_1(X) \le \lambda_2(X) \le \cdots \le \lambda_n(X) \le +2
$$
\item The smallest eigenvalue is $0$, that is, $\lambda_1(X) = 0$. It realized by the eigenvector
$f_1(v) = \sqrt{d_v}$.
\item In the case that $X$ is $k$-regular, so $\Delta = I - (1/k) A$, the eigenvalues of
the adjacency matrix $A$ are $\hat{\lambda}_j(X) = k - k \ \lambda_j(X), j= 1 , \cdots, n$ with
$$
-k \le  \hat{\lambda}_n(X) \le \hat{\lambda}_{n-1}(X) \le \cdots \le \hat{\lambda}_2(X) \le \hat{\lambda}_1(X) = k
$$
\
\item For a subset $A$ of vertices of $X$ with complement $B = V(X) - A$,  define the volume of $A$ to be
$$
vol(A) = \Sigma_{v  \in  A} \ \ d_v
$$
and the edge set of $A$, $\partial A$ to be the set of directed edges with initial point in $A$ and
terminal point in $B$, That is,
$$
\partial A = \{ \vec{e} \ | \  i(\vec{e})  \in A \ and \ t( \vec{e}) \in B \ \}
$$
Then the Cheeger type estimate states:
$$
\lambda_2(X) \le |\partial A| \ ( \frac{1}{vol(A)} + \frac{1}{vol(B)})
$$
In particular, for $vol(A) \le vol(B)$, one has
$$
\lambda_2(X) \le 2 \frac{|\partial A|}{vol(A)}
$$
This is  a geometrically based isoperimetric inequality where $|\partial A|$ serves
as the measure of the ``surface area'' of the boundary of $A$ and $vol(A)$ serves
as the volume.
\end{itemize}

\vspace{.3in}

The definition of the transmission graph Laplacian is modeled on the above
explicit approach.
For each vertex $v $ of $X$ let $n(v)$ be a non-negative integer and
form the vector space $C^{n(v)}$ of column vectors of length $n(v)$.
This integer $n(v)$ is called the rank of the vertex $v$.
The vectors $a \in C^{n(v)}$ are to be thought of as  a choice of data
residing at the
vertex $C^{n(v)}$.

 Concretely, for vertices $v,w$ of the graph $X$,  a linear mapping from $C^{n(v)}$ to   $ C^{n(w)}$  is given
 by  a $n(w)$ by $ n(v)$ matrix with complex entries.
A transmission system is  a function  from the directed edges of the graph $X$
to matrices
$$
P :  E^+(X)  \rightarrow  \ matrices
$$
where the matrix $P(\vec{e}$ is to be a $n(w)$ by $  n(v)$ matrix giving the
linear mapping  of $C^{n(v)}$ to  $C^{n(w}$ by $ a \mapsto P(\vec{e}) a$.
[That is, the $i^{th}$ entry in this column is $ \Sigma{1 \le j \le n(v)} \ P(\vec{e})_{i,j} \ a_j $.]

$P(\vec{e})$ are to  be thought of as a means of transmitting  data at the initial point
$v = i(\vec{e})$ recorded as a column vector  $ a $  in $ C^{n(v)}$ to data at the terminal point
$w = t(\vec{e})$ by sending $ a \mapsto P(\vec{e}) \ a$ which lies in $C^{n(w)}$

Consequently, for each pair of vertices, say $v_i,v_j$, there is an associated linear mapping of $C^{n(v_i)} $ to
$C^{n(v_j)} $
obtained by summing, $P(\vec{e})$, over directed edges $\vec{e}$  with initial point $v_i$ and terminal point $v_j$. It is
given by adding the $n(v_j)$ by $  n(v_i)$ matrices over these directed edges.
$$
A^P_{i,j}  =   \Sigma_{ \vec{e} \ with \ i(\vec{e}) = v_i, \ t(\vec{e}) = v_j   } \  P(\vec{e}) \ .
$$
In the case that $N=1$ and $P(\vec{e}) = Id$ for all directed edges, this linear mapping is just multiplication
by the adjacency matrix element $A_{i,j}$.

Now grouped together these $n$ by $n$ blocks  form an $N$ by $N$ matrix with $ N = \Sigma_{v \in V(X)} v(v)$
$$
A^P = \{ \ A^P_{i,j} \ | \ 1 \le i,j \le n \ \}
$$
which is called the adjacency matrix of the transmission system $P$.

Correspondingly, define the normalized transmission matrix $N(A)^P$ by grouping together the blocks
$$
N(A)^P_{i,j} = \frac{1}{ \sqrt{d_{v_i} \ d_{v_j} } } \  \Sigma_{ \vec{e} \ with \ i(\vec{e}) = v_i, \ t(\vec{e}) = v_j } \  P(\vec{e})
$$

In these terms, the transmission graph Laplacian is defined by the $nN$ by $nN$ matrix
$$
\Delta^P = Id - N(A)^P  = Id - \{ \ N(A)_{i,j}^P \ \}
$$

Alternatively, a more elegant way to describe $\Delta^P$ is as follows. Let $C^0(X,C^N)$ denote
the set of mappings of $V(X) $ to $C^N$,
$$
C^0(X,C^N) = Map( V(X) , C^N)
$$
This is a complex vector space under addition and scalar multiplication by elements of $C$.
 It is  isomorphic to the complex vector space $C^{nN} = C^N \oplus C^N \oplus \cdots \oplus C^N$ via the bijection
$$
f \mapsto   (f(v_1), f(v_2), \cdots, f(v_n))
$$
For $f : V(X) \rightarrow C^N$ define $(\Delta^P f) \in C^0(X, C^N)$ via
$$
(\Delta^P f)(w) = f(w) - \Sigma_{ \vec{e} \ with \ t(\vec{e} ) = w } \
\  \frac{1}{\sqrt{d_w \ d_{i(\vec{e} )} } } \ P(\vec{e}) f(i(\vec{e}))
$$
Under the identification $C^0(X,C^N) \cong C^{nN}$ this agrees with the above explicit description of $\Delta^P$.

\vspace{.3in}
Simple estimates in section \ref{sectestimates1} prove the theorem:

\begin{thm}    \label{thmrange}
a) If $\lambda $ is an eigenvalue of $\Delta^P$, then the real and imaginary parts of
$\lambda$ are bounded by
$$
\begin{array}{l}
| \Re( \lambda) - 1 | \le | \lambda - 1| \le max_{\ \vec{e} \in E^+(X) \ }  \ \ ||P(\vec{e}) || \\
\\
\Im( \lambda) \le (1/2) \  max_{\ \vec{e} \in E^+(X) \ }  \ \ ||P(\vec{e})^\star - P(op(\vec{e}))  || \\
\end{array}
$$
In particular, if $P$ is a (Hermitian) symmetric transmission system, that is,
$ P(op(\vec{e}))= P(\vec{e})^\star$ for all directed edges $\vec{e}$, then
$\Delta^P$ is self-adjoint and has real eigenvectors.

In particular, if $P$ is a strictly unitary transmission system, that
is  (Hermitian) symmetric and invertible, the eigenvalues of $\Delta^P$ are all real
and lie in the range $[0,+2]$.
\end{thm}

Let $P$ be a (Hermitian) Symmetric transmission system, Since the eigenvalues are real
they may be listed in increasing order with appropriate multiplicities. Let them be
in increasing order
$$
\lambda_1(X,P) \le \lambda_2(X,P) \le \cdots \le \lambda_{nN}(X,P)
$$

Let $f_1 : V(X) \rightarrow C^N$ be a choice of non-zero eigenvector for the lowest real eigenvalue
$\lambda_1(X,P)$ of $\Delta^P$, that is,
$$
\Delta^P \ f_1 = \lambda_1(X,P) \ f_1 \ .
$$
[In the classical case $f_1(v) = \sqrt{d_v}$ with $\lambda_1((X) = 0$.]

For a subset $A$ of vertices of the graph $X$ with complement $B = V(X) - A$, define the volume of $A$ by the formula
$$
vol_P(A) = \Sigma_{ v \in A} \ \  ||f_1(v)||^2
$$
[In the classical case this is $\Sigma_{v \in A} \ d_v = vol(A)$.]

Recall that the edge set of $A$, $\partial A$, consists of directed edge
s $\vec{e}$ with
initial point in $A$ and terminal point in $B$. Define
$$
|\partial A|_P = \Sigma_{ \ \vec{e} \in \partial A  } \  \ \frac{|f(i(\vec{e})) |  }{ \sqrt{ d_{  i(\vec{e} ) } }  }
\ \frac{|f(t(\vec{e})) |  }{ \sqrt{ d_{  t(\vec{e} ) } }  }
   $$
   [In the classical case, this is exactly $|\partial A|$.]

   In these terms the generalization of the upper Cheeger type estimate is:

\begin{thm}   \label{thmCheegerupper1}
For a (Hermitian) symmetric transmission system $P$, the difference of second minus  first smallest eigenvalues,
$\lambda_2(X,P)- \lambda_1(X,P)$, is bounded above by
$$
\lambda_2(X,P)- \lambda_1(X,P) \le K \ |\partial_P(A)| | ( \frac{1}{vol_P(A)} + \frac{1}{vol_P(B)} )
$$
where $K = max_{ \vec{e} \ } \ ||P(\vec{e})||$ and $B= V(X)-A$ is the complement of $A$.

In particular, for $vol_P(A) \le vol_P(B)$, this difference is bounded by
$ 2  K \ \frac{ |\partial_P(A)|}{vol_P(A)} $.
\end{thm}

A more elegant formulation in terms of capacities is given in section \ref{sectCapacities}.
As seen from the proof in section \ref{sectestimates1} this upper bound can be
strengthened.

Now replacing the values of $\frac{f(v)}{\sqrt{d_v}}$ by their maximums and minimums this theorem
implies the weaker version:

\begin{thm}   \label{thmCheegerupper2}
For a (Hermitian) symmetric transmission system $P$, the difference of second minus  first smallest eigenvalues,
$\lambda_2(X,P)- \lambda_1(X,P)$, is bounded above by
$$
\lambda_2(X,P)- \lambda_1(X,P) \le  \ |\partial(A)| | ( \frac{1}{vol(A)} + \frac{1}{vol(B)} ) K L
$$
where $K = max_{ \vec{e} \ } \ ||P(\vec{e})||$ and $B= V(X)-A$ is the complement of $A$
and
$$
L = [max_{ v \in V(X) \ } \ \frac{ |f(v)|}{\sqrt{d_v}} ]^2/ [ min_{ v \in V(X) \ with \ |f(v)| \neq 0 } \ \frac{ |f(v)|}{\sqrt{d_v}}]^2
$$
\end{thm}

As seen,  these have precisely the form of an isoperimetric inequality as desired.

\vspace{.3in}

An important case of transmission graph Laplacians
occurs in studying random walks. Here there are
$n$ states, say $\{v_1,v_2,\cdots, v_n\}$.
At each ordered pair of states, say $v_i, v_j$, a probability
$Prob(v_i,v_j)$ is assigned as the probability that state $v_i$
will transmutate into state $v_j$. The probabilites are to be real
non-negative and to satisfy the constraint $\Sigma_{j} \ Prob(v_i,v_j) = +1$,
a conservation of probability condition. Let $\Gamma$ be the transition
probability matrix with $j,i$ entry $Prob(v_i,v_j)$. Then the $j,i$
entry of $\Gamma^N$ is the probability that a random walk of length
$N$ will pass from state $v_i$ to state $v_j$.

Now let $G$ be the graph with vertices $\{v_i\}$ and edges
$\{ v_i,v_j\}$ if either of $Prob(v_i,v_j)$ or $Prob(v_j,v_i)$ are
non-zero. Recall that here the transmission system is defined
in terms of directed edges where one sums over incoming directed
edges in forming the transmission graph Laplacian.
Define a transmission system $P$ assigning as above to a directed
edge $\vec{v_iv_j}$ the probability $Prob(v_i,v_j)$ if $v_i \neq v_j$.
If  $Prob(v,v) ) \neq 0$ for a vertex $v$, set the probability of transversing the
loop in either direction to be $(1/2)Prob(v,v)$ so the sum over all directed
edges $\Sigma_{\vec{v} \ with \ i(\vec{v}) = w} \ P(\vec{v}) +1$ as desired.
 Let $T$ be the degree mapping which
sends the vertex $v$ to its degree $d_v$. Let $Q$ be the transmission system
which sends the directed edge $\vec{v}$ to $d_{t(\vec{v})} \ P(\vec{v})$.

With this notation, the set of eigenvalues of the transition matrix $\Gamma$
is computed equivalently by:
$$
\begin{array}{l}
eigen(\ \Gamma\ ) = eigen(\ T^{1/2} \ \Gamma \ T^{-1/2} \ ) = 1 - eigen( \ Id - T^{1/2} \ \Gamma \ T^{-1/2}\ ) \\
= 1 - eigen( T^{-1/2} \ ( Id - T \ \Gamma)\  T^{-1/2} \ )
\end{array}
$$
That is, subtracting each of the eigenvalues of $\Gamma$ from one, one gets precisely
the eigenvalues of the transmission system $Q$. The system $Q$ is (Hermitian) symmetric
precisely when
$$
d_{v_j} \ Prob(v_i,v_j) = d_{v_i} \ Prob( v_j,v_i)
$$
Such probabilistic systems are called reflexive in the literature \cite{Chung1}.

\section{A Symmetric Graph Transmission System for Associations.} \label{sectassociation}

Let $X$ be a simple graph with no loop edges. Let $V(x)$ denote the
set of vertices. Write $ a \sim_X b$ if the vertices $a,b$ are
adjacent in the graph $X$.

An \textbf{association}  $A$ is a mapping of $V(X)$ to the power set $2^{V(X)}$.
That is, $A$ gives for each vertex $v$ of $X$ a subset $A(v) \subset V(X)$.
We regard the vertices, say $w$, of $A(v)$ as the associates of $v$.
Write $a \Rightarrow_A w$ in this instance, that is, if $ w \in A(a)$.

To a graph $X$ and an association $A$, define a symmetric transmission system
and thus a graph Laplacian in the following manner.

\vspace{.3in}

For a subset,say $S$ of $V(X)$, let $F(S)$ be the vector space of functions $f
: S \rightarrow R$ of the set $S$ to the real numbers. For two
subsets,say $S,T$, there is a natural linear mapping
$$
{\cal{I}}(S,T) : \ F(S) \rightarrow F(T)
$$
defined as follows: Send $f: S \rightarrow R$ to the function, say
$\hat{f} : T \rightarrow R$, given by $\hat{f}(b) = 0$ if $b \in
T-S$ and $\hat{f}(b) =f(b)$ if $b \in S \cap T$.

Note that ${\cal{I}}(S,T)$ is defined by passing the information
$f(b)$ for vertices $b$ common to both $S$ and $T$ from the
function $f$ to the function $\hat{f}$.

Endowing $F(X)$ with the real inner product, $<f,g> = \Sigma_{x \in
X} \ f(x) g(x)$ by restriction the subspaces $F(S),F(T)$ inherit
real inner products denoted by $<.,.>_S, <.,.>_T$ respectively. It
is seen that the two linear mappings ${\cal{I}}(S,T),
{\cal{I}}(T,S)$ are adjoint. That is, for $f \in F(S), g \in F(T)$
we have
$$
< {\cal{I}}(S, \ T) f, g>_T = < f, {\cal{I}}(T, \ S) g>_S \ .
$$

\vspace{.3in}
 With these notations in mind,  introduce  the
 a  real symmetric transmission system on the
 graph $X$  given  an association $A$ on the set of vertices $V(X)$
 as follows:

For each vertex, say $a $ of $X$,  associate the vector space
$$
G[a] := F(A(a)),
$$ and to each directed edge $\vec{ab}$ of $X$
 associated the natural linear mapping
$$
{\cal{I}}(A(a)   , \ A(b))  \ : \ G[a] \rightarrow G[b]
$$
of course here $ a \sim_X b $.

This is a real symmetric transmission system, by the above observations about
${\cal{I}}(S,T)$, so has an associated
transmission graph Laplacian $\Delta$ with real eigenvalues.
This Laplacian $\Delta$ is easily described directly in the following
manner.

By definition $\Delta$ acts linearly on the $0-cochains $ in
this setting is just the direct sum over the vertices of the vector
spaces $G[a]$:
$$
C^0 \ =  \ \bigoplus_{a \in X} \ G[a]
$$
Such an element is just a list $\{K(a)\}$ with $K(a) \in G[a]$.
These $K(a)$ are functions $K(a) : A(a) \rightarrow R$.

\vspace{.3in} For vertices $a,v$ let $N(a,v)$ be the integer
$$
N(a,v)\ =  \ \#  \{ b \ |   \ a \sim_X b \ and \ b \Rightarrow_A v \}
$$
This is the number of vertices $b$ which are adjacent to $a$ in the graph $X$
while $v$ is an associate of the vertex $b$, i.e., $v \in A(b)$.
Note $N(a,v) \ge 1$  if there is such a $b$.

 In these terms the linear mapping  $\Delta$ sends $\{K(a)| a \in X
\}$ to $\{L(b)| b\in X \}$ precisely when for all $b \in X$ and $z
\in A(b)$

$$
L(b)(z) = H(b)(z) - \Sigma_{a\in X \ with\  a \sim_X b \ and \
a \Rightarrow_A z } \ \frac{1}{\sqrt{N(a,z) \ N(b,z)}} \ K(a)(z)
$$
where if the set summed over is vacuous the sum is set equal to $0$
by convention.
Note that if the set $\{a\in X \ with\  a \sim_X b \ and \ a \Rightarrow_A z \}$ is non-empty, then $ z \in A(a) \cap A(b)$ and $a
\sim_X b$, so $N(a,z) \ge 1$ and  $N(b,z) \ge 1$ as desired in forming this definition.

An elementary analysis shows that $\Delta$ has $0$ as an eigenvalue
with explicitly specified eigensolution,  that
all the eigenvalues lie in the range $[0,2]$, and a Cheeger type
estimate for the smallest non-zero eigenvalue, called here $\lambda(X,A)$.
The proofs are carried out in the rest of this section.

\vspace{.3in}
By the above noted symmetry of ${\cal{I}}(S, \ T)$, the linear mapping
$\Delta$ is self-adjoint and so has real eigenvalues.

Define for each vertex $a$ the element $\Phi(a) \in G[a] = F(A(a))$
given by the function $\Phi(a) : A(a) \rightarrow R$,
$$
\Phi(a)(v) = \sqrt{N(a,v)} \ for \  a \Rightarrow_A v
$$

It is claimed that the 0-cochain $ a \mapsto \Phi(a)$ is in the
kernel of $\Delta$.

To see this directly  compute $\Delta(\{ \Phi(a)\})$ at a point, say $b$. Now for $
w \in A(b)$,
$$
\begin{array}{l}
(Id - \Delta) \Phi \ (b)(w) \\
= \Sigma_{a \ with \ a \Rightarrow_A w \ and \ a \sim_X b} \
   \frac{1}{\sqrt{N(a,w) \ N(b,w)}} \Phi(a)(w) \\
   =\Sigma_{a \ with \  a \Rightarrow_A w \ and \ a \sim_X b } \
   \frac{1}{\sqrt{N(a,w) \ N(b,w)}} \  \sqrt{N(a,w)} \\
   = \frac{1}{\sqrt{N(b,w)}} \ \Sigma_{a \ with \  a \Rightarrow_A w \ and \ a \sim_X
   b} \ 1 \\
   = \sqrt{N(b,w)} = \Phi(b)(w) = Id \ \Phi(b)(w).
   \end{array}
   $$
In the exceptional case that the sets $\{ a \ with \  a \Rightarrow_A w
\ and \ a \sim_X
   b\}$ are vacuous for all $a$, then each term is by convention
   zero and moreover $N(b,w)$ is zero so the formula still obtains.

\vspace{.3in}

 To show that the eigenvalues of $\Delta$ lie in the range $[0,2]$,
 since they are real it suffices to show that the norm of $Id -
\Delta$ is at most one.

Let $\{K(a)\}$ be a non-vanishing  $0$-cochain, so in particular,
\newline  $K(a) : A(a) \rightarrow R$ for each vertex $a$. We estimate
from above $< \{K(a)\} ,(Id- \Delta)( \ \{K(a)\})> $ exactly as in Chung's book \cite{Chung1},
using the elementary inequality $rs \le (1/2)( r^2 + s^2) $ for real
numbers $r,s$:
$$
\begin{array}{l}
< \{K(a)\} ,(Id- \Delta) \ \{K(a)\}> = \Sigma_b \ <K(b), [Id - \Delta] (\{K(a)\})(b) > \\
=\Sigma_b \ \Sigma_{(a,w)  \ with \ w \in A(a), w \in A(b),
\ and \ a \sim_X
   b} \  < K(b)(w), \  \ \frac{1}{\sqrt{N(a,z) \ N(b,z)}} \ K(a)(w) >\\
   \le (1/2) \  \Sigma_b \ \Sigma_{(a,w) \ with \  w \in A(a), w \in A(b), \ and \ a \sim_X
   b} \     [ \frac{|K(b)(w)|^2}{N(a,z)}  + \frac{|K(a)(w)|^2}{N(b,z)} ]\\
   = (1/2) \  \Sigma_b \ \Sigma_{(a,w) \ with \  w \in A(a), w \in A(b), \ and \ a \sim_X
   b}\ \frac{|K(b)(w)|^2}{N(a,z)} \\
   + (1/2)  \ \Sigma_a \ \Sigma_{(b,w) \ with \  w \in A(a), w \in A(b), \ and \ b \sim_X
   a}\ \frac{|K(a)(w)|^2}{N(b,z)}   \\
   =(1/2) \ \Sigma_{(b,w) \ with\  w \in A(b) } \ |K(b)(w)|^2 \\
   + (1/2) \ \Sigma_{(a,w) \ with \ w \in A(a) }  \ |K(a)(w)|^2 \\ = \Sigma_{(a,w) \ with \ w \in A(a)}  \
   |K(a)(w)|^2 = |\{K(a)\}|^2
\end{array}
$$

\vspace{.3in}
 In a similar manner, one  gets a Cheeger type of estimate
giving an upper bound for the next non-negative eigenvalue
$\lambda_1(X,A)$ for any two disjoint sets, $A,B$, of vertices.

\vspace{.3in}

Given disjoint sets $A,B$ of vertices.

Now form the 0-cochain, $\{\Psi(a)\}$ defined by
$$
\begin{array}{l}
\Psi(a) = -S \ \Phi(a) \ for \ a \in A ,\\
\Psi(b) = +R \ \Phi(b) \ for \ b \in B \\
\Psi(x) = 0 \ for \ x \notin A \cup B
\end{array}
$$
where $R,S$ are chosen as
$$
\begin{array}{l}
R =\Sigma_{a \in A} \ |\Phi(a)|^2 =  \Sigma_{a\in A, v \in A(a)} \ N(a,v)   \\
S  =\Sigma_{b \in B} \ |\Phi(b)|^2 = \Sigma_{b \in B, w \in A(a)} \ N(b,w)
\end{array}
$$

With these definitions, one obtains,
$$
\begin{array}{l}
<\{\Phi(x)\}, \{ \Psi(y)\}> \\
= \Sigma_{a\in A} \ (-S) |\Phi(a)|^2 + \Sigma_{b \in B} \ (R) |\Phi(b)|^2 \\
= -RS + RS = 0.
\end{array}
$$

By the mini-max principle, since $\Delta$ is real self-adjoint, the next largest eigenvalue,
say $\lambda_1(X,A)$, of $\Delta$, is then estimated from above by the Ritz-Rayleigh quotient:
$$
 \lambda_1(X,A) \le \frac{< \{\Psi(x)\}, \Delta \ \{\Psi(x)\}>}{ < \{\Psi(x)\}, \{\Psi(x)\}>}
 $$

 One obtains simply,
 $$
 \begin{array}{l}
 < \{\Psi(x)\}, \{\Psi(x)\}> \\
 = \Sigma_{a \in A}\ S^2 |\Phi(a)|^2 +  \Sigma_{b \in B}\ R^2 |\Phi(b)|^2 \\
 = R S^2 + S R^2 = RS(R+S).
 \end{array}
 $$

Similarly, using the formula $\Delta ( \{\Phi(a)\}  ) = 0$, one obtains,
$$
\begin{array}{l}
< \{\Psi(x)\}, \Delta(\{\Psi(x)\}) > \\
= \Sigma_{a \in A} \ < (-S)\Phi(a) , \ \Delta(\{\Psi(x)\})(a)> \\
+ \Sigma_{b \in B} \ < (R) \Phi(b) , \ \Delta(\{\Psi(x)\})(b) > \\
= -\Sigma_{a \in A} \Sigma_{ b \in B \ with \ a\sim_X b, v \in A(a)\cap A(b) } \ (-S)(R-(-S))\
 \frac{< \Phi(a)(v), \Phi(b)(v)>}{\sqrt{N(a,v) N(b,v)}} \\
 - \Sigma_{b \in B} \Sigma_{ a \in A \ with \ b\sim_X a, v \in A(a)\cap A(b) } \ (R)(-S-R)\
 \frac{< \Phi(b)(v), \Phi(a)(v)>}{\sqrt{N(b,v) N(a,v)}} \\
 = S(R+S) \ \Sigma_{a\in A, b \in B \ with \  a\sim_X b, v \in A(a)\cap A(b)} \ 1 \\
 +R(R+S) \  \Sigma_{a\in A, b \in B \ with \ a\sim_X b, v \in A(a)\cap A(b)} \ 1 \\
 = (R+S)^2 \ \Sigma_{a\in A, b \in B \ with \ a\sim_X b, v \in A(a)\cap A(b)} \ 1 .
 \end{array}
 $$

 Hence, one gets  a Cheeger type estimate:
 $$
 \begin{array}{l}
 \lambda_1(X,A) \le \frac{< \{\Psi(x)\}, \Delta \ \{\Psi(x)\}>}{ < \{\Psi(x)\}, \{\Psi(x)\}>} \\
  =  2 \ ( \frac{1}{R} + \frac{1}{S})\ \Sigma_{a\in A, b \in B \ with \ a\sim_X b, v \in A(a)\cap A(b)} \ 1
  \end{array}
 $$

\vspace{.3in}
Note that the quantity $\Sigma_{a\in A, b \in B \ with \ a\sim_X b, v \in A(a)\cap A(b)} \ 1$ is precisely
the sum over edges $\vec{ab}$ with $a \in A$, $b \in B$ with $a$ adjacent to $b$ of the number of
vertices $v$ which are associates of \textbf{both} $a$ and $b$. Thus, the proper way to measure
the ``cohesion'' of the graph with respect to the association $A$ is to form the above quantity,
$2 \ ( \frac{1}{R} + \frac{1}{S})\ \Sigma_{a\in A, b \in B \ with \ a\sim_X b, v \in A(a)\cap A(b)} \ 1$.
Once done this way, it provides an upper estimate for the eigenvalue $\lambda_1(X,A)$.
That is, the above prescription provides the correct way of meaningfully measuring the ``cohesion'',
for example in social networks with associations.

\vspace{.3in}
Additionally, on can introduce for each vertex a density or weighting function $f(v): A(v) \rightarrow R$
on the set associated to the vertex $v$ and incorporate this into the above definiitions. Again similar
results hold. This increased flexibility is useful in various  applications.

\section{Quantum Mechanical Example of Transmission Systems: Invertible and (Hermitian) Symmetric:} \label{sectquantum}

Here an example of a (Hermitian) symmetric transmission system
is given which is not also strictly unitary.
These arise naturally in quantum mechanical settings, most simply in the  one dimensional
transmission of left and right going free particles
through a potential well.  This section  provides an
 invertible transmission system, $F(\vec{e})$, which when multiplied by
a fixed matrix $J$ gives a new transmission system
$J \ F(\vec{e})$ which is (Hermitian) symmetric. Hence,
this modified system will have real eigenvalues.
This provides physically natural   examples of (Hermitian) symmetric
transmission systems beyond the obvious strictly unitary ones.

The idea is to consider  for any directed edge, say $\vec{PQ}$,
 the propagation of particles, say of spin $n/2$,
in a one dimensional manner along this edge. Suppose moreover,
there is a potential barrier between $P$ and $Q$ while our spin $n/2$
particles move freely left and right nearby $P$, $Q$, as  the
barrier is absent there. In this setting quantum mechanics gives a transmission
matrix relating the flux of the free particles at $P$ to the
flux of free particles at $Q$. It is this matrix that yields
 an invertible transmission system,  for our graph as follows:

More explicitly, let $x$ be a coordinate of the edge with $P,Q$
at positions, say $x(P), x(Q)$ with $x(P) < x(Q)$. Then for
$A,B \in C^{n+1}$ by $n+1 = 2(n/2)+1$ for fixed energy $E$
there are the free fluxes
$$
\Psi_L(x) = e^{ ik (x-x(P))} \  A + e^{-ik(x-x(P)x)}\ B  \ with\ A, B \in C^{n+1} \ for \ x \ near\ by \ x(P)
$$
nearby $P$ of spin $n/2$ particles. One imagines that for fixed column vector
$(A,B)^t$  in $ C^{n+1} \oplus C^{n+1} = C^{2(n+1)}$ these free fluxes appear in each out-going edge from $P$
in a neighborhood of  $P$.
Similarly, at $Q$ there is the data $(C,D)^t$ with associated wave
$$
\Psi_R(x) = e^{ ik (x-x(Q))} \  C + e^{-ik(x-x(Q))} \  D  \ with \ C, D \in C^{n+1} \ for \ x \ near\ by \ x(Q)
$$

Consider  free particles of spin $n/2$ moving in one dimension towards and away from
a potential barrier located in the region $x(P) < x < x(Q)$ and continuing
through to the region $x>Q$. In these two regions one has the
representation for fixed energy $E $ particles  of right and left moving fluxes:
$$
\begin{array}{l}
\Psi_L(x) = e^{ ik (x-x(P))} \  A + e^{-ik(x-x(P))}\ B  \ with\ A, B \in C^{n+1} \ for \ x <P \\
 \Psi_R(x) = e^{ ik (x-x(Q))} \  C + e^{-ik(x-x(Q))} \  D  \ with \ C, D \in C^{n+1} \ for \ x >Q
 \end{array}
 $$
 Here we normalize the left, right going waves to have phase 0 at $P,Q$  in the
 regions $x<P, x>Q$ respectively.

 By quantum mechanics, the arbitrary values of $(A,B)^t$ determine uniquely the values
 $(C,D)^t$ and the relation is recorded by multiplication by the  transmission matrix, say $M$, with
$$
\left( \begin{array}{c} C \\ D \end{array} \right)
= M \ \left( \begin{array}{cc} A \\ B \end{array} \right),  \ with \  M = \left( \begin{array}{cc} R & S \\ U & V \end{array} \right)
$$
where $M$ is a $2(n+1) $ by $2(n+1)$ invertible matrix since spin $n/2$ particles take their
values in $C^{2(n/2)+1} = C^{n+1}$ with its standard inner product.
Here the inverse $M^{-1}$ reverses this process sending $(C,D)^t$ to $(A,B)^t$.
The associated ordinary differential
equations are   of second order, whence the free solutions nearby $P$ uniquely extend to
solutions of the ODE throughout the interval, in particular in a neighborhood of $Q$.

\{ Actually, in standard texts \cite{Merzbacher} the conventional transmission matrix is $M^{-1}$,
going from $(C,D)^t$ to $(A,B)^t$,
but this is immaterial.\}

Now conservation of probability, a general feature of quantum mechanics,
demands that $|A|^2 - |B|^2 = |C|^2 - |D|^2$   or equivalently
via $|A|^2 - |B|^2  = <(A^\star,B^\star) , I_{n+1,n+1} (A,B)^t> $ for
$$
I_{n+1,n+1} = \left( \begin{array}{cc} 1 & 0 \\
                              0 & -1 \end{array} \right)
                              $$
                              Consequently, the transmission matrix
                              $M$ must preserve the form $|A|^2 - |B|^2$,
                              so it is necessarily lies in the Lie group $U(n+1,n+1)$
or equivalently, $M$ satisfies the added constraint
$$
I_{n+1,n+1} = M^\star \ I_{n+1,n+1} \ M
$$

\{ It is of interest to note that the commutator $[I\ ,\ M ] = I \ M \ I^{-1} \ M^{-1}
= I \ M \ I M^{-1} = I \ M \ M^\star \ I = (I \ M) \ (I \ M)^\star$. In particular,
the commutator $[I , M]$
$$
[ I \ , \ M]= (I \ M) \ (I \ M)^\star
$$
is self-adjoint, so has real eigenvalues. Here $I = I_{n+1,n+1}$.
In physical measurements, one expects to easily recover
the element $M \in U(n+1,n+1)$ only upto left and right multiplication
by $U(n+1) \times U(n+1)$. It should be noted that the eigenvalues
of the commutator $[I,M]$ factor thought the double
coset space $U(n+1) \times U(n+1)\backslash U(n+1,n+1) / U(n+1) \times U(n+1)$
as desired so  measurement of its real eigenvalues may be envisioned
in practice.\}

Now if one represents  $P,Q$ as an edge of a graph $X$, the  data
$(A,B)$ located at $P$ as yields   on  each out going edge, say $\vec{PQ}$,
 the specified free  wave $\Psi_L(x)$ nearby $P$.  This wave    is transmitted forward and reflected back
by a potential well between $P$ an $Q$  and once propagated through the barrier
emerges as a free
left and right moving current as $\Psi_L$, then
it is natural to associate to this directed edge $\vec{PQ}$
the transmission matrix which records the relation
of out-going at P, in-coming at P mapping to
out-going at Q, in-coming at Q. Since right going at
P = out-going at P while right going at Q = in-coming at
Q, this matrix is precisely:
$$
F( \vec{P Q}) =  Perm \ M  \ with \ Perm = \left( \begin{array}{cc} 0 & 1 \\ 1 & 0 \end{array} \right)
$$
for the permutation matrix $Perm$.

Working from $Q$ to $P$, for outgoing at $Q$ to outgoing at $P$ gives
the matrix $M^{-1} \  Perm =(Perm \ M)^{-1}$ and so
$$
F(\vec{QP}) = M^{-1} \ Perm = F(\vec{PQ})^{-1} ,
$$
This association $\vec{e} \mapsto F(\vec{e})$ defined for an directed edge
$\vec{e}$
 is an invertible transmission system for the given graph $X$ recording
 its transmission along its  directed edges $\{ \vec{e} \}$.

Now define for each  directed edge $\vec{e}$,
$$
G(\vec{e}) =  i \ I_{n+1,n+1}  \  \ F(\vec{e}) \ with \ i = \sqrt{-1}
$$

It is claimed that $\vec{e} \mapsto G(\vec{e})$ is a (Hermitian) symmetric
transmission system.

This is precisely the computation, utilizing $ I \ M^{-1} = M^\star \ I$,
$$
\begin{array}{l}
G(\vec{QP})=
[  i \ I \  \ F(\vec{QP}) ] \\
 = [ i  \ I \ \ M^{-1} Perm  ] =  [ i  \  M^\star \ I \ Perm ]  \\
         = [ -i \ Perm \ I \  M  ] ^\star = [ i \ I \ Perm \ M]^\star \\
         = [ i \ I \ F(\vec{PQ})]^\star = G(\vec{PQ})^\star
         \end{array}
         $$

  A natural idea is to use this transmission of incoming and outgoing particles to record a
transmission matrix as determined explicitly by general quantum mechanical principles.

This example motivates the ``quantum transmission graph Laplacian'' defined as
the transmission graph Laplacian for this (Hermitian) symmetric transmission
system $G(\vec{e})$. This operator is self-adjoint and has real eigenvalues. Presumably,
in discussions of transmitting
qubits this circle of ideas has relevance.

\section{The Dual of a Graph and its Transmission Systems.} \label{sectdual}

For a graph $X$ with vertex set $V(X)$ and edge set $E(V)$, the dual of
$X$ is a graph, $X^\star$ with one vertex, say $D(e)$ for each undirected edge
and a edge from $D(e)$ to $D(f)$ whenever there is a vertex of $X$
appearing as an end point of both the edges $e,f$.

Consequently, a transmission system on the dual $X^\star$ consists
of an assignment of vector spaces, say $W(e)$ to each edge of $X$
and for each pair of edges, $e,f$ with common vertex of a linear
mapping $W(e) \rightarrow W(f)$ and visa versa. It is natural
to think of this in the following way. A vector in $W(e)$ records
the data being transmitted along the edge $e$ and the linear
mapping $W(e) \rightarrow W(f)$ gives the rule for transmission
of this data through the common vertex, or node.

To be more explicit, suppose that $X$ is a simple graph without loops
and whose edges are determined by their endpoints.
 Let $X$ have $m$ vertices $v_1,\cdots, v_m$ and
$n$ edges $e_1,\cdots, e_n$. Form the incidence matrix, $Incidence$, a $m$ by $n$ matrix,
with entry $+1$ at the $i,j$ entry if the vertex $v_i$ is a end point of the
edge $e_j$. Consequently,  the $m$ by $m$ symmetric matrix $(Incidence) \cdot (Incidence)^t $
has $i,j$ entry equal to $d(v_i)$ the degree of $i$ if $i=j$ and equal to
$+1$ if there is an edge $e$ with end points $v_i,v_j$ if $i \neq j$ and zero otherwise.
That is,
$$
(Incidence) \cdot (Incidence)^t = T + Adjacency(X)
$$
where $Adjacency(X)$ is the adjacency matrix of $X$.

Similarly, the $n$ by $n$ matrix $(Incidence)^t \cdot (Incidence)$
has diagonal entries all $+2$ and the off diagonal terms
equal to the adjacency matrix of the dual graph $X^\star$.
$$
(Incidence)^t \cdot (Incidence) =  2 \ Id + Adjacency(X^\star)
$$

Note that if $X$ is $k$-regular, then $X^\star$ is $2(k-1)$ regular
In particular, since the non-zero eigenvalues of $(Incidence) \cdot (Incidence)^t $
and $(Incidence) \cdot (Incidence)^t$ coincide,  in this $k$ regular case, one
gets easily that
$$
k \ \lambda_2(X) = 2(k-1) \ \lambda_2(X^\star)
$$
for the classical normalized graph Laplacian, here
$Id - (1/k) Adjacency(X), Id - (1/(2(k-1))\ Adjacency(X^\star)$
respectively.

\section{Estimates 1: Range of Values and Cheeger's Upper Bound Generalized} \label{sectestimates1}

Proof of Theorem \ref{thmrange}:

Introduce the Hermitian inner product on $C^0(X,C^N)$ given by
$$
<f,g> = \Sigma_v \   < f(v),g(v)>
$$
for $f,g : V(X) \rightarrow C^N$. Here $<a,b>$ for $a,b \in C^N$ is the standard Hermitian inner
product on $C^N$.

Let $f$ be an eigenvector for an eigenvalue $\lambda$ of the transmission graph Laplacian
$\Delta^P$ with $f$ chosen so that $||f||^2 = <f,f> = 1$. Then
using the inequality $|rs| \le (1/2) ( |r|^2 + |s|^2)$, one obtains
$$
\begin{array}{l}
|\lambda -1| =  |< f, \Delta^P f  - f>| \\
= | - \Sigma_{\vec{e} \ with \ i(\vec{e}) = v, t(\vec{e}) =w \ }
\frac{1}{\sqrt{d_v \ d_w}} \ <f(w), P(\vec{e}) \ f(v) > |
\\
\le  \Sigma_{\vec{e} \ with \ i(\vec{e}) = v, t(\vec{e}) =w\ }
\frac{1}{\sqrt{d_v \ d_w}} \ |<f(w), P(\vec{e}) \ f(v) >| \\
\le  \Sigma_{\vec{e} \ with \ i(\vec{e}) = v, t(\vec{e}) =w\ }
\frac{1}{\sqrt{d_v \ d_w}} \  ||P(\vec{e})|| \ ||f(w)||\  ||f(v) || \\
\\
\le (1/2)  \Sigma_{\vec{e} \ with \ i(\vec{e}) = v, t(\vec{e}) =w\ }
 \ || P(\vec{e})|| \ ( \frac{||f(v)||^2}{d_v} + \frac{||f(w)||^2}{d_w}  ) \\
\\
\le max_{ \vec{e}\ } \ ||P(\vec{e}) || \ ( \Sigma_v \ ||f(v)||^2 )= max_{ \vec{e}\ } \ ||P(\vec{e}) ||
\end{array}
$$

Also, taking conjugates,
$$
\begin{array}{l}
\overline{\lambda} - 1 = - \Sigma_{\vec{e} \ with \ i(\vec{e}) = v, t(\vec{e}) =w \ }
\frac{1}{\sqrt{d_v \ d_w}} \ \overline{<f(w), P(\vec{e}) \ f(v) > }\\
= - \Sigma_{\vec{e} \ with \ i(\vec{e}) = v, t(\vec{e}) =w \ }
\frac{1}{\sqrt{d_v \ d_w}} \ \overline{<P(\vec{e})^\star f(w), \ f(v) > }\\
= - \Sigma_{\vec{e} \ with \ i(\vec{e}) = v, t(\vec{e}) =w \ }
\frac{1}{\sqrt{d_v \ d_w}} \ < f(v), P(\vec{e})^\star f(w) > \\
while \\
\lambda - 1 = - \Sigma_{\vec{e} \ with \ i(\vec{e}) = v, t(\vec{e}) =w \ }
\frac{1}{\sqrt{d_v \ d_w}} \ <f(v), P(op(\vec{e})) \ f(w) > \\
\end{array}
$$
Subtracting gives
$$
 2i \ \Im(\lambda) = + \Sigma_{\vec{e} \ with \ i(\vec{e}) = v, t(\vec{e}) =w \ }
\frac{1}{\sqrt{d_v \ d_w}} \ <f(v), [P(\vec{e})^\star - P(op(\vec{e}))] \ f(w) >$$
 by the
above method this yields the desired result \newline $  \Im(\lambda)
\le (1/2) \ max_{ \vec{e} \ } \ ||P(\vec{e})^\star - P(op(\vec{e})) ||$.

\vspace{.3in}
Proof of Theorem \ref{thmCheegerupper1}:

Let $f_1$ be an eigenvector for the lowest real eigenvalue $\lambda_1(X,P)$
of the transmission graph Laplacian, not necessarily normalized.
Define a function $ g : V(X) \rightarrow C^N$ by
$$
\begin{array}{l}
g(v) = S \ f_1(v) \ if \ v \in A \\
g(v) = - R \ f_1(v) \ if \ v \in B = V(X) -A \\
\end{array}
$$
where $R =\Sigma_{v \in A} \ ||f_1(v)||^2$ and $S = \Sigma_{v \in B} \ ||f_1(v)||^2$.

With these choices, $<f_1,g> = \Sigma_v \ < f_1(v),g(v)>
= \Sigma_{v \in A} \ S \ ||f_1(v)||^2 - \Sigma_{v \in B} \ R \ ||f_1(v)||^2 =SR-RS=0$.
Consequently, $g$ is orthogonal to the  eigenvector $f_1$ with the smallest eigenvalue
$\lambda_1(X,P)$ of $\Delta^P$.
Since $\Delta^P$ is self-adjoint for (Hermitian) symmetric transmission system,
the Raleigh-Ritz principle applies to estimating the second smallest eigenvalue
$\lambda_2(X,P)$ giving:
$$
\lambda_2(X,P) \le \frac{ < g , \Delta^P g>}{ < g, g>}
$$

Here $<g,g> = \Sigma_{v \in A} S^2 \ ||f_1(v)||^2 + \Sigma_{v \in B} \ R^2 \ ||f_1(v)||^2
= S^2 \ R  + R^2 \ S $. This gives
$$
\begin{array}{l}
 \lambda_2(X,P) - \lambda_1(X,P) \le  \frac{ < g , \Delta^P g> - \lambda_1(X,P)  \ (S^2\ R + R^2\ S)}{<g,g>} \\
 = \frac{[ \Sigma_{v\in A} < g(v) , (\Delta^P g)(v)>  - <S f_1(v), S \ (\Delta^P f_1)(v)> ]}{<g,g>} \\
 +  \frac{[ \Sigma_{v\in B} < g(v) , (\Delta^P g)(v)>  - <R f_1(v),R \ (\Delta^P f_1)(v)]>  }{<g,g>}
 \end{array}
 $$

Now for $v \in A$, $ < g(v) , (\Delta^P g)(v)>  - <S f_1(v), S \ (\Delta^P f_1)(v)> \newline
= <S f_1(v), (\Delta^P g)(v) - S \ (\Delta^P f_1)(v)>$. Here by $ v\in A$,
 $$
 \begin{array}{l}
 (\Delta^P g)(v)
= S f_1(v) - \Sigma_{\vec{e} \ with \ t(\vec{e}) = v, i(\vec{e}) = w, \ with \ w \in A} \ \frac{1}{\sqrt{d_v \ d_w}} \ S \ P(\vec{e}) f_1(w)\\
 \hspace{.3in} + \Sigma_{\vec{e} \ with \ t(\vec{e}) = v, i(\vec{e}) = w, \ with \ w \in B} \ \frac{1}{\sqrt{d_v \ d_w} } \
  R \ P(\vec{e}) f_1(w)\\
 S(\Delta^P f_1)(v)
= S f_1(v) - \Sigma_{\vec{e} \ with \ t(\vec{e}) = v, i(\vec{e}) = w, \ with \ w \in A} \  \frac{1}{\sqrt{d_v \ d_w}}  \
 S\  P(\vec{e}) f_1(w)\\
 \vspace{.3in} - \Sigma_{\vec{e} \ with \ t(\vec{e}) = v, i(\vec{e}) = w, \ with \ w \in B} \  \frac{1}{\sqrt{d_v \ d_w}} \
 S \ P(\vec{e}) f_1(w)\\
 \end{array}
 $$

 Consequently, the contributions from edges connecting points internal to $A$ vanish, giving for $ v \in A$,
 $ (\Delta^P g)(v) - S(\Delta^P f_1)(v) \newline  = (R+S) \ \Sigma_{\vec{e} \ with \ t(\vec{e}) = v, i(\vec{e}) = w, \ with \ w \in B }
 \ P(\vec{e})\ f_1(w)$. Similarly, for $w\in B$, \newline $ (\Delta^P g)(w) - S(\Delta^P f_1)(w)
  =(R+S) \  \Sigma_{\vec{e} \ with \ t(\vec{e}) = w, i(\vec{e}) = v, \ with \ v \in A}
 \ P(\vec{e})\ f_1(v)$. In toto, these yield, summing over $A$ and $B = V(X)-A$,

 $$
 \begin{array}{l}
 [\lambda_2(X,P) - \lambda_1(X,P)] \ <g,g>  \\
 \le  \ S(R+S) \ \Sigma_{\vec{e} \ with \ i(\vec{e}) = v \in A , i(\vec{e}) = w \in B}
     \ \frac{1}{\sqrt{d_v \ d_w}} \ <f_1(v), P(\vec{e}) f_1(w)> \\
    + R(R+S)  \ \Sigma_{\vec{e} \ with \ i(\vec{e}) = v \in A , t(\vec{e}) = w \in B}
      \ \frac{1}{\sqrt{d_v \ d_w}} \ <f_1(w), P(\vec{e}) f_1(v)> \\
 = (R+S)^2  \ \Sigma_{\vec{e} \ with \ t(\vec{e}) = v \in A , i(\vec{e}) = w \in B}
    \  \ \frac{1}{\sqrt{d_v \ d_w}} \ <f_1(v), P(\vec{e}) f_1(w)>
    \end{array}
    $$
        where the (Hermitian) symmetry was used to rewrite
    $ <f_1(v), P(\vec{e}) f_1(w)> = <P(\vec{e})^* f_1(v), f_1(w)>
    = <P(op(\vec{e}) f_1(v), f_1(w)>
    = \overline{ < f_1(w), P(op(\vec{e}) f_1(v)>}$ and real parts are taken on both sides to collect the
    terms.

By standard estimates, this last gives the desired inequality:
$$
\begin{array}{l}
[\lambda_2(X,P) - \lambda_1(X,P)] \\
\le \frac{(R+S)^2}{RS(R+S)} \Sigma_{\vec{e} \ with \ t(\vec{e}) = v \in A , i(\vec{e}) = w \in B} \
\frac{1}{\sqrt{d_v \ d_w}} \ |<f_1(v), P(\vec{e}) f_1(w)>| \\
\le \frac{(R+S)^2}{RS(R+S)} \Sigma_{\vec{e} \ with \ t(\vec{e}) = v \in A , i(\vec{e}) = w \in B}
\ ||P(\vec{e})|| \frac{||f_1(v)||}{\sqrt{d_v}} \ \frac{||f_1(w)||}{\sqrt{d_w}} \\
\le [ max_{\ \vec{e} \ } \ ||P(\vec{e})|| ] \ \Sigma_{\vec{e} \ with \ t(\vec{e}) = v \in A , i(\vec{e}) = w \in B} \
\frac{||f_1(v)||}{\sqrt{d_v}} \ \frac{||f_1(w)||}{\sqrt{d_w}} \  \ (\frac{1}{R} + \frac{1}{S})
\end{array}
$$

The proof of Theorem \ref{thmCheegerupper2} follows directly from the above.
The first inequality is a strengthened form of this result.

\section{Estimates 2: Eigenvalues and Diameters} \label{sectNilli}

Let vertices, $v,w$, be called adjacent if there is a directed
edge in $X$ from $v$ to $w$. In this case write $ v \sim w$.
 The distance between two vertices $v,w$
is $L= d(v,w) $ where $v_1,v_2,\cdots, v_L$ is the shortest list of vertices
with $v=v_1, w= v_L$,  and $v_i \sim v_{i+1}$ for $ 1 \le i \le L-1$.
The diameter of $X$  is the maximum over all such distances.
For a vertex $v$ of $X$, let $\hat{d}(v)$ be the number of vertices
of $X$ \textbf{distinct} from $v$ which are adjacent to $v$. That is,
$\hat{d}(v)$ does not count the possible loops at the point $v$.

A graph is called $k$-regular if $d(v) = k$ for each vertex $v$.

For a $k$-regular graph, Nilli bounded the second
eigenvalue of the classical normalized Laplacian
as follows in terms of diameters:

\begin{thm}  \label{thmNilli1} \cite{Nilli1}
Let X be a $k$-regular graph.
 If the diameter of X is $\ge 2 b + 2 \ge 4$,
then the second eigenvalue $\lambda_2(X)$ is bounded above as follows:
$$
\lambda_2(X) \le 1 - (1/k) [ 2 \sqrt{k-1} - \frac{2 \sqrt{k-1} -1}{b}]
$$
\end{thm}

For $X$ a $k$-regular graph this is equivalently: the  adjacency matrix $A$ has second to the highest eigenvalue,
that is , $k - k \ \lambda_2(X)$, is bounded below  by
$$
k - k \ \lambda_2(X) \ge [ 2
\sqrt{k-1} - \frac{2 \sqrt{k-1} -1}{b}]
$$

A $k$-regular graph is called
a Ramanujan graph if
$$
k - k \lambda_2(X)  \le  2 \ \sqrt{k-1}
$$
that is, if its adjacency matrix $A$ has second to the highest eigenvalue,
 $k - k \ \lambda_2(X)$, is
at most $2 \ \sqrt{k-1}$. These are graphs of significant interest in applications,
for example, to robust codes and switching theory \cite{Chung1, Chung2,Gnutzman, Lubotzky, Murty, Sarnak, Sipser}.

It is to be noted that for the classical normalized graph Laplacian
the matrix obtained from a graph $X$  and a graph obtained from $X$
by deleting any loops are identical. This is not the case
for our general transmission Laplacian since the the transmission
data for loops may not be the identity; however, the Cheeger
type constants are insensative to these loops.

\vspace{.3in}
For a graph $X$ and a vertex $v $ of $X$, let $\hat{d}(v)$ be the number
of vertices adjacent to $v$  but not equal to it. That is, $\hat{d}$
is the degree $d_v$ minus twice the number of loops at $v$.
The analogue in the transmission setting of Theorem \ref{thmNilli1} is:

\begin{thm}  \label{thmNilli2}
Let X be a graph with $\hat{d}(v)  \le k$ for all vertices
and  diameter of X  $\ge 2 b + 2 \ge 4$. Let
$P$ be a (Hermitian) symmetric transmission system on $X$.
Let $q = (k-1)$, and let $Ave(d)$ be the weighted average
of the degrees $d_v$
$$
Ave(d) = \Sigma_{v} \ d_v \ (||f_1(v)||^2/d_v )/ \Sigma_v \ (||f_1(v)||^2/d_v).
$$
Then the difference of the first and second eigenvalues of
the transmission graph Laplacian satisfies:
$$
\lambda_2(X,P) - \lambda_1(X,P) \le (1/Ave(d)) \
   (K/M) \ ( 1 + qM^2 - 2 \sqrt{q}M  + \frac{ 2 \sqrt{q} M -1}{b} )
$$
where $K
 = max_{\vec{e}} \ ||P(\vec{e})||$ and \newline
$M = max_{\vec{e} \ with f_1(i(\vec{e})) \neq 0 } \ \ (||f(t(\vec{e}))||/\sqrt{d_{t(\vec{e})} })/
(||f(i(\vec{e}))||/\sqrt{d_{i(\vec{e})} } )$.

\end{thm}

Note that in the classical case $f_1(v) =d_v$ so $K=M=1$ and for $X$
$k$-regular, $Ave(d) = k$ and  the Nilli result is recovered

\vspace{.3in}
Proof:
Following Murty \cite{Murty}, chose $u,v \in X$ be such that $d(u,v)  \ge 2b+2$.
Define sets for $ i \ge 0$ by
$$
\begin{array}{l}
U_i =  \{x \in V(X) \ | \  d(x,u) = i \} \\
V_i =  \{x \in V(X) \ | \  d(x,v) = i \}.
\end{array}
$$
The sets $U_0,U_1,\cdots, U_b, V_0,V_1, \cdots, V_b$ are disjoint,
for otherwise by the triangle inequality one gets $d(v,w) \le 2b$, which is
a contradiction.
Moreover, no vertex of $U = \bigcup \ U _{i=0}^b$ is adjacent to a vertex of
$V = \bigcup_{j=0}^b \ V_j$ since otherwise $d(u,v) \le 2b+1$ which is again
a contradiction. For each vertex in $U_i$ at least one edge leads to  a vertex of $U_{i-1}$
and at most $q = k-1$ to a vertex of  $U_{i+1}$ ( for $  i \ge 1$).

Let $f_1$ be the eigenvector for the first eigenvalue $\lambda_1(X)$ and $g_1 = T^{-1/2} \  f_1$,
that is, $g_1(v) = d_v^{-1/2}  \ f_1(v)$. Since the normalized graph Laplacian is of the
form $ T^{-1/2} \ ( T - A) \ T^{-1/2}$, $\Delta \ f_1 = \lambda_1(X) \ f_1$  implies that $ \lambda_1(x) \ Tg_1 = (T-A) g_1$
where $A$ is the generalized adjacency matrix for the transmission graph system.

Thus, for $ i \ge 1$
$$
\begin{array}{l}
  \Sigma_{v \in U_{i+1}} \ ||g_1(v)||^2 \le qM^2 \  \ \Sigma_{v \in U_i} \ ||g_1(v)||^2 \\
and \ similarly \\
  \Sigma_{v \in V_{i+1}} \ ||g_1(v)||^2 \le qM^2 \  \ \Sigma_{v \in V_i} \ ||g_1(v)||^2.
  \end{array}
  $$

  Let $r = M \sqrt{q}$. Let $g$ be the function defined by $g(v) = F_i \ g_1(v)$ for
  $ v \in U_i$ and $g(w) = G_i \ g_1(w)$ for $w \in V_i$ for $ 0 \le i \le b$
and zero otherwise. Here $F_i,G_i$
  are specified by $F_0 = \alpha, G_0 = \beta, F_i = \alpha r^{-(i-1)}, G_i = \beta r^{-(i-1)}$
  for $ 1 \le i \le b$. Let $g= T^{-1/2} f$ define $f$.

 Choose $\alpha, \beta $ so that $<f_1,f> = 0$, or equivalently $<g_1, T g> = 0$. Note $F_0=F_1,G_0=G_1$.

  By definition $<g,g> = A+ B$ with $A= \Sigma_{i=0}^b \ ||F_i||^2 \ ( \Sigma_{v \in U_i} \ ||g_1(v)||^2)$
  and $B= \Sigma_{j=0}^b \ ||G_j||^2 \ ( \Sigma_{w \in V_i} \ ||g_1(w)||^2 )$. By  the mini-max principle \newline
  $ (<g, (T-A)g> -
  \lambda_1(X) <g,Tg>)/<g,Tg> = ( <f, \Delta f> - \lambda_1 < f, f>)/<f,f> \ \ge \  \lambda_2 - \lambda_1$.

  Again $<g, (T-A) \ g> - \lambda_1(X)  < g, Tg> = C + D$ with \newline
  $C = \Sigma_{i=0}^b \ F_i ( \Sigma_{v \in U_i} \ < g_1(v), ((T-A) g)(v) - F_i \ ((T-A) g_1)(v) > $
  and \newline
  $D = \Sigma_{i=0}^b \ G_i ( \Sigma_{w \in V_i} \ < g_1(w), ((T-A) g)(w) - G_i \ ((T-A) g_1)(w) > $.

  Since $ ((T-A) h)(v) = h(v) - \Sigma_{\vec{e} \ with \ t(\vec{v}) = v} \ P(\vec{e}) \  h(i(\vec{e}))$, for
  $v \in U_i$
$$
\begin{array}{l}
F_i \  < g_1(v), ((T-A) g)(v) - F_i \ ((T-A) g_1)(v) > \\
= \Sigma_{\vec{e} \ with
     \ t(\vec{e}) = v \ and \ i(\vec{e}) \in U_{i+1} } \ \ F_i (F_{i+1}-F_i)\ < g_1(i(\vec{e})), P(\vec{e}) g_1(t(\vec{e}))>\\
   +  \Sigma_{\vec{e} \ with
     \ t(\vec{e}) = v \ and \ i(\vec{e}) \in U_{i-1} } \ \ F_i (F_{i-1}-F_i) \ < g_1(i(\vec{e})), P(\vec{e}) g_1(t(\vec{e}))>
     \end{array}
     $$
     assuming $ 1 \le i \le b$. Since $F_0 = F_1$ the corresponding  sum vanishes for $v \in U_0$.
  For an directed edge $\vec{e}$ with initial vertex in $U_i$ and terminal vertex in $U_{i+1}$,
  the contribution of this edge to $C$ comes in two parts one from the directed edge
  $\vec{e}$, that is, $F_i (F_{i+1}-F_i)\ < g_1(i(\vec{e})), P(\vec{e}) g_1(t(\vec{e})))>$, and the
  other from the oppositely directed edge, that is $F_{i+1} (F_{i}-F_{i+1})\ < g_1(i(\vec{e})), P(op(\vec{e}))
 g_1(t(\vec{e})))> = F_{i+1} (F_{i}-F_{i+1}) \ \overline{ < g_1(t(\vec{e})), P(\vec{e}) g_1(i(\vec{e})))>}$ by Hermitian
symmetry. The sum of these two contributes to the real part as $(F_i -F_{i+1})^2 \Re{< g_1(i(\vec{e})), P(\vec{e}) g_1(t(\vec{e})) >}$.
Since $C$ is real, such pairs contribute to $C = \Re{C}$ by  a magnitude at most $(F_i -F_{i+1})^2 \ K \ M ||f(t(\vec{e}))||^2$.
Since there are at most $q$ edges with initial vertex in  $U_{i+1}$ which terminate in $U_i$, the sum
over these edges contributes at most $ q\ K M \ (F_i -F_{i+1})^2 ( \Sigma_{v \in U_i} \||g_1(v)||^2)$
to and upper estimate of $|C|$. When $i=b$ by convention $F_{b+1} = 0$.
 Similar considerations for $D$ apply yielding the upper estimates:
$$
\begin{array}{l}
|C| \le \Sigma_{i=1}^{b-1} \ q KM \ (F_i-F_{i+1})^2  \ (\Sigma_{v \in U_i} \ ||g_1(v))||^2 ) \\
 \hspace{.5in} + \ q KM \ F_b^2  \ (\Sigma_{v \in U_b} \ ||g_1(v)||^2 ) \\
 = (q KM) \ ( r^{-1} -1)^2 ( \Sigma_{i=1}^b \ r^{-2(i-1)}  \ (\Sigma_{v \in U_i} \ ||g_1(v)||^2)) \  \alpha^2\\
       + (q KM) \ (-r^{-2}+2r^{-1}) r^{-2(b-1)} \ (\Sigma_{v \in U_b} \ ||g_1(v)||^2)) \  \alpha^2 \\
       = (K/M) \ ( r -1)^2 ( \Sigma_{i=1}^b \ r^{-2(i-1)}  \ (\Sigma_{v \in U_i} \ ||g_1(v)||^2)) \  \alpha^2\\
       + (K/M) \ (1-2r) r^{-2(b-1)} \ (\Sigma_{v \in U_b} \ ||g_1(v)||^2)) \  \alpha^2 \\
 |D| \le \Sigma_{j=1}^{b-1} \ q KM \ (G_i-G_{i+1})^2  \ (\Sigma_{w\in V_j} \ ||g_1(w))||^2 ) \\
 \hspace{.5in} + \ q KM \ G_b^2  \ (\Sigma_{w \in V_b} \ ||g_1(w)||^2 ) \\
 = (K/M) \ ( r -1)^2 ( \Sigma_{i=1}^b \ r^{-2(i-1)}  \ (\Sigma_{v \in U_i} \ ||g_1(v)||^2)) \  \alpha^2\\
       + (K/M) \ (2 r -1) r^{-2(b-1)} \ (\Sigma_{v \in U_b} \ ||g_1(v)||^2) ) \ \alpha^2 \\
\end{array}
$$

  Since $r^{2(i-1)} (\Sigma_{v \in U_i} \ ||g_1(v)||^2) \ge r^{2i} (\Sigma_{v \in U_{i+1}} \ ||g_1(v)||^2)$
  for $ b \ge i \ge 1$ and similarly for $V_j$,
  one has
 $$
  \begin{array}{l}
  |C| \le   (K/M) \ (r-1)^2 \ ( A - \alpha^2) + (K/M)( 2 r-1) ( A - \alpha^2)/b \\
   \le (K/M) \ ( 1 + qM^2 - 2 \sqrt{q}M  + \frac{ 2 \sqrt{q} M -1}{b}) A \\
  |D| \le   (K/M) \ (r-1)^2 \ ( B - \beta^2) + (K/M)( 2 r-1) ( B - \beta^2)/b \\
  \le (K/M) \ ( 1 + qM^2 - 2 \sqrt{q}M  + \frac{ 2 \sqrt{q} M -1}{b}) B
  \end{array}
  $$

Consequently, $\lambda_2- \lambda_1 \le \frac{C+D}{\Sigma_v \ d_v ||g_1(v)||^2} = (1/Ave(d)) \ \frac{C+D}{A+B} \newline
\le \
(1/Ave(d)) \  (K/M) \ ( 1 + qM^2 - 2 \sqrt{q}M  + \frac{ 2 \sqrt{q} M -1}{b})$
as claimed.

\section{Reinterpretation in terms of Capacities. Cheeger Constants.}  \label{sectCapacities}

Given a set say $S$, of vertices, for each directed edge, $\vec{v_i,v_j}$, possibly with $v_i = v_j$
assign a real non-negative number, $C(v_i,v_j)$, a capacity. For each ordered pair of
vertices, say $v_i,v_j$ assign a weight $D(v_i.v_j)$, real non-negative numbers, not all zero.
 Then for any subset $A \subset S$
one may assign the quotient
$$
F(A) = \frac{  \Sigma_{\vec{v_i,v_j}}  \  C(\vec{v_iv_j}) \ | \  1_A ( v_i) - 1_A(v_j) \ |}
         { \Sigma_{i,j} \ D(v_i,v_j)\ | \  1_A ( v_i) - 1_A(v_j) \ | }
         $$
where $1_A$ is the characteristic function of the set $A$ which vanishes outside $A$ and
is identically one inside. This quotient $F(A)$ is a measure of the amount of capacity
of the edges which pass from inside $A$ to outside $A$ or visa versa divided by the
sum of the weights of these pairs. Minimizing over the subsets $A$ yields a measure
of the capacity of the total system.

Now in the case $X$ is a graph and $D(\vec{v})$ is zero or one according to
$\vec{v}$ being a directed edge of $X$ while $E(v_i,v_j)$ is identically one,
$$
F(A) = \frac{\partial A}{|A| \ |V(X)-A|}
$$
so the minimum,  once  multiplied by $|V(X)|$, is a variant of  the classical Cheeger constant.
[In the literature the Cheeger constant is defined by minimizing $F(A)\ |V(X)-A| = \frac{\partial A}{|A|}$
over subsets $A$ of size at most $(1/2)|V(X)|$. But if  $|A| \le (1/2) |V(X)|$
one has
$$
\frac{1}{|A|} \le \frac{|V(X)|}{|A| \ |V(X)-A|} \le \frac{2}{|A|}
$$
so minimizing over these $A$ and minimizing $F(A) \ |V(X)|$ over all $A$
are essentially equivalent.

By looking at the proof of $\ref{thmCheegerupper1}$, one sees that the upper
estimate for the difference in eigenvalues $\lambda_2(X,P)-\lambda_1(X,P)$
is the minimum of $F(A) \|V(X)|$ for the capacities and weights
$$
\begin{array}{l}
C(\vec{v}) = \frac{  | <\ f_1( t(\vec{v}) )  ,\  P(\vec{v}) \ f_1(i(\vec{v}))  \ >  }{\sqrt{ d_{i(\vec{v})} \ d_{t(\vec{v})} } } \\
D(i(\vec{v}), t(\vec{v}) )=  d_{i(\vec{v})} \  d_{t(\vec{v}) }
\end{array}
$$

These capacities may be replaced by the weaker version
$$
C'(\vec{v}) = K \frac{ f_1(i(\vec{v})}{\sqrt{ d_{i(\vec{v})} }  }\ \frac{f_1(t(\vec{v})}{\sqrt{d_{t(\vec{v})}} }
$$
as in the statement of theorem \ref{thmCheegerupper2}.

\section{Cayley graphs: Their Eigenvalues.} \label{sectCayley}

A graph is $k$-regular if each vertex has degree $k$.
There is a simple procedure for constructing $k$-regular graphs using group
theory. This can be described as follows. Let $G$ be a finite group and $S$ a $k$-
element multi-set of $G$. That is, $S$ has $k$ elements where one allows repetitions.
Suppose that $S$ is symmetric in the sense that $s \in S $ implies $s^{-1} \in S$  (with
the same multiplicity). Now construct the graph $(X,G; S)$  by having the vertex
set to be the elements of $G $ with  $(x,y,s)$ a directed edge from $x$ to $y$  if and only if $s= x^{-1} y$
lies in $S$. Let the same edge with the opposite direction be defined to be $(y,x,s^{-1})$ which naturally satisfies
the desired constraint $s^{-1} = y^{-1} x \in S$.
 Since
$S$ is allowed to be a multi-set, $(X,G; S)$ may have multiple edges.
If G is abelian, the eigenvalues of the Cayley graph are easily determined
as follows.

\begin{thm} ( adapted from \cite{Murty})  \label{thmCayley} Let G be a finite abelian group and S a symmetric subset of $ G $
of size $k$. Then the eigenvalues of the graph Laplacian of the $k$-regular graph $(X,G; S)$ are given by
$$
1- (1/k)\ [ \Sigma_{ s \in S} \ \chi(s)]
$$
as $\chi$  ranges over all the irreducible characters of $G$.
\end{thm}

This theorem allows explicit control of all eigenvalues of these special
Caley graphs. In particular, it allows one to search for  Ramanujan graphs. Those are $k$-regular
graphs with  second eigenvalue $\lambda_2$ satisfying :
$$
\lambda_2 \le   1 - (1/k)\ [ 2 \sqrt{k-1} ]
$$
which is close to the optimal upper bound for  $\lambda_2$.
See Theorem \ref{thmNilli1}.

\vspace{.3in}
Now suppose that  a function $F$ associates to each element, say $s$, of $S$
a square $N$ by $N$ invertible matrix, $F(s)$, with the added property:
$$
F(s^{-1}) = F(s)^\star
$$
Then one may define a (Hermitian) symmetric transmission system on the
Cayley graph $G$ by defining $P(x,y,s) = F( s) = F(x^{-1} y)$. This will be (Hermitian) symmetric
transmission system.

A  transmission analogue of the above theorem is:

\begin{thm}  \label{thmCayley} For $G$ an abelian group and $S$, $F$ as above,
 the eigenvalues of the transmission graph Laplacian of the $k$-regular graph $(X,G; S)$
 with transmission system $P$  are given by
$$
1 - (1/k) [ \hat{\lambda(\chi, F)}]
$$
where $\lambda(\chi,F)$ ranges over  the $N$ eigenvalues of the matrix
$\Sigma_{s \in S} \ \chi(s) \ F(s) $ and $\chi$ over the irreducible characters
of $G$.
\end{thm}

\vspace{.3in}
Proof: Let $\delta_S(q)$ equal to zero if $q \notin S$ and equal to the multiplicity
of $q$ in $S$ if $ q \in S$. Since the graph is $k$-regular, the eigenvalues of the
graph Laplacian $\Delta^P$ are equal to $1- (1/k) \hat{\lambda}$ where $\hat{lambda}$
are the eigenvalues of the transmission graph adjacency matrix $A^P$ of the system $P$.
As usual set $C^)(X,C^N) = Map(G,C^N)$, in this  notation the $A^P :  \rightarrow
C^0(X,C^N)$ is given by
$$
(A^P \ f)(x) = \Sigma_{g \in G} \ \delta_S(x^{-1} y) \ F(x^{-1}g) \ f(g)
$$

Now let $e_i, i = 1, \cdots, N$ be the standard basis for  $C^N$ and for each character
$\chi :G \rightarrow C$ define the element $v(\chi,i) \in C^0(X,C^N)$ by
$$
v(\chi,i)(g) = \chi(g) \ e_i
$$
As $\chi$ ranges over the irreducible characters of $G$ and $i$ ranges from  $ 1$ to $ N$
these form a basis for $C^0(G,C^N)$ by $G$ abelian.

Now the effect of $A^P$ on this basis element is, replacing $x^{-1} g $ by $s$
$$
\begin{array}{l}
A^P( v(\chi,i)) = \Sigma{b \in G} \ \delta_S(x^{-1} y) \ F(x^{-1}g) \ \chi(g) \ e_i \\
         = \chi(x) \ [ \Sigma_{x \in S} \ F(s) \chi(s) \ e_i] \\
         = [\Sigma_{s \in S} \ F(s) \chi(s)] v(\chi ,i)
         \end{array}
         $$
The theorem now follows.

\section{Eigenvalues and Graph Collapses and Amalgamations.} \label{sectcollapse}

A useful technique in graph theory is to collapse a graph, that is, simplify it
by identifying vertices. This procedure has excellent properties with respect to
the eigenvalues of the graph Laplacian and as we will see, to its generalization, the transmission graph
Laplacian.

Here the collapse of a graph $X$ with transmission data $P$ is defined very easily.
The collapsed graph, called  $X'$,  has precisely the same directed edges as the un-collapsed
graph $X$ with the same transmission system, but the vertices of $X$, $ V(X)$
are identified to get a smaller set, say $V'$. That is, there is a surjection
$\Psi : V(X) \rightarrow V'$ and two vertices, say $v,w$ are identified if
$\Psi(v) = \Psi(w)$. Let $X'$ denote the collapsed graph.
The condition for this to be sensible is that if two vertices, say $v,w$ are to
be identified, then they must have equal bit ranks $n(v) = n(w)$.
This approach to collapsing is
easier and more transparent than that usually one
\cite{Chung1}.

\begin{thm} \label{thmcollapse}
If $X'$ is obtained by collapsing $X$ and the transmission system for $X$
is (Hermitian) symmetric, then the induced transmission  system for $X'$ is also
(Hermitian) symmetric and
the eigenvalues of $X'$, say
 $ \mu_1 \le \mu_2 \le \mu_3 \cdots $, are related to those of
$X$,  say
$\lambda_1 \le \lambda_2 \le \lambda_3 \le \cdots $,
[both counted with multiplicities by the inequalities] by the inequalities
$$
 \lambda_j \le \mu_j
$$
\end{thm}

Proof: The surjection $\Psi$ induces an inclusion of vector spaces:
$$
\begin{array}{l}
I : C^0(X',C^N) = Map(V',C^N) \rightarrow C^0(X,C^N) = Map(V(X) ,C^N) \\
 f \mapsto  \{ v \mapsto f(\Psi(v)) \\
 \end{array}
 $$
 Note that the degree of a vertex $v \in X'$, say $\overline{d}_v$, is given
 in terms of the degrees of the pre-images in $X$ by:
 $$
 \overline{d}_v =  \Sigma_{w \in V(X) \ with \ \Psi(w)=v \ } \ d_w
 $$
 In particular, the map $I$ has the property that for maps
 $ f, g : V' \rightarrow C^N$, one has
 $$
 \Sigma_{v \in V'} \ \overline{d}_v < f(v), g(v) > = \Sigma_{w\in V(X) \  } \ d_w < f(\Psi (w)), g(\Psi (w))>,
 $$
 so $I$ is an isometry for these inner products.

 Let $ T : C^0(X,C^N) \rightarrow C^0(X,C^N)$ be the grading mapping sending $f$ to
 $\{ w \mapsto  d_w \ f(w)\} $ and $T' : C^0(X',C^N) \rightarrow C^0(X',C^N) $ sending $f$
 by $\{ v \mapsto \overline{d}_w \ f(w) \}$.

 In these terms the transmission graph Laplacian $\Delta^P$ is of the form
 $$
 \begin{array}{l}
 \Delta^P = T^{-1/2} ( T - L) T^{-1/2}  \\
 where \\
 (Lf)(v) = \Sigma_{\vec{e} \ with \ t(\vec{e}) = v \ } \  P(\vec{e}) \ f(i(\vec{e}) )
 \end{array}
 $$

 Now for $g = T^{-1/2} f$ for $ f : V(X) \rightarrow C^N$, the Rayleigh-Ritz
 quotient for the operator $\Delta^P$ may be reexpressed in term of
 $L$ since
 $$
 \begin{array}{l}
 <f, \Delta^P f> = <f, T^{-1/2} (T - L) T^{_1/2} f>=  < g, (T-L) g> \\
 and \ <f,f>  = \Sigma_{w \in V(X)\ } \ d_w ||g(w)||^2
 \end{array}
 $$
 as
 $$
 \frac{<f, \Delta^P f > }{<f,f>} = \frac{  < g, (T-L) g>}{ \Sigma_{w \in V(X) \ } \ d_w ||g(w)||^2 }
 $$
 Similarly, for $X'$ the formula under the translation $g = T'^{-1/2} f$ is
 $$
 \frac{<f, \Delta'^P f > }{<f,f>} = \frac{  < g, (T'-L') g>}{ \Sigma_{w \in X' \ } \ \overline{d}_w ||g(w)||^2 }
 $$

 It is claimed that the transmission graph Laplacians,
 $\Delta^p, \Delta'^P$ for $X,X'$ respectively, satisfy the
 compatibility condition for $ g : V' \rightarrow C^N$
 $$
 [\star] \hspace{.3in} < (g \cdot \Psi), (T-L)  \ (g \cdot \Psi)>  = < g, (T'-L') g>
 $$

 Granting this,   the Rayleigh-Ritz quotients for $X'$ may be
 computed equivalently in the subspace $I(C^0(X',C^N)) \subset C^0(X,C^N)$
 using  $\frac{  < g, (T'-L') g>}{ \Sigma_{w \in X' \ } \ \overline{d}_w ||g(w)||^2 }$,
 consequently the appropriate mini-max estimates on $X'$ are greater than
 those of $X$ which are taken over the complete space $C^0(X,C^N)$. This proves
 the inequalities
 $$
 \lambda_j \le \mu_j
$$

To prove $[\star]$ is suffices to do so for an elementary collapse in which just
two distinct vertices, say $P,Q$ are identified. Let $V(X) - \{ P,Q \} = B$ and
the vertices of $X'$ be denoted by $[v]$ for $  v \in B$ and the remaining vertex, say $X$, which
is the result of identifying $P$ and $Q$.

Chose $ g : V' \rightarrow C^N$.
Now part of the left hand side of $[\star]$ is the sum
$<  (g \cdot \Psi)(P), ((T-L)  \ (g \cdot \Psi))(P)> + <  (g \cdot \Psi)(Q), ((T-L)  \ (g \cdot \Psi))(Q)>
= < g(X),((T-L)  \ (g \cdot \Psi))(P) + ((T-L)  \ (g \cdot \Psi))(Q) > $. Here
$$
\begin{array}{l}
((T-L)  \ (g \cdot \Psi))(P) + ((T-L)  \ (g \cdot \Psi))(Q) \\
= (d_P + d_Q)  \ g(X) - \Sigma_{\vec{e} \ with \ t(\vec{e}) \in \{ P,Q\}, i(\vec{e}) \notin \{P,Q \}  \ } \  P(\vec{e})\ g( i(\vec{e})) \\
- \Sigma_{\vec{e} \ with \ t(\vec{e}) \in \{P,Q\}, i(\vec{e}) \in \{P,Q\} \ } \ P(\vec{e}) \ g(X) \\
= \overline{d}_X \ g(X) -  \Sigma_{\vec{e} \ with \ \Psi(t(\vec{e}))  = X \ }  \ P(\vec{e}) g(\Psi(i(\vec{e}))) \\
= (\Delta'^P g)(X)
\end{array}
$$
Together this yields $<  (g \cdot \Psi)(P), ((T-L)  \ (g \cdot \Psi))(P)> + <  (g \cdot \Psi)(Q), ((T-L)  \ (g \cdot \Psi))(Q)>
= < g(X), (\Delta'^P g)(X)>$.

In a similar manner, if $v \notin \{P,Q\}$, then
$$
\begin{array}{l}
<  (g \cdot \Psi)(v), ((T-L)  \ (g \cdot \Psi))(v)> \\
= < g(v), d_v \ g(v) - \Sigma_{\vec{e} \ with \ t(\vec{e}) = v , i(\vec{e}) \notin \{P,Q\} \ } \ P(\vec{e}) \ g(i(\vec{e}) ) \\
  -  \Sigma_{\vec{e} \ with \ t(\vec{e}) = v , i(\vec{e}) \in  \{P,Q\} \ } \ P(\vec{e}) \ g(X)  \\
  = < g(v) , (\Delta'^P g)( v) >
  \end{array}
  $$
Adding these together gives the required formula $[\star]$.

\vspace{.3in}
\textbf{Amalgamation of graphs and transmission systems.}

Given a graph $X$ one can amalgamate the edges of the graph by making a new graph $X'$ with
the same vertices but with exactly one  edge between two vertices, say
$v,w$ not necessarily distinct whenever   there is a directed edge from
$v$ to $w$, or equivalently visa versa. Now if $P$ is a transmission
system for $X$, let the amalgamated system $P'$ be defined as follows:
If $v,w$ are distinct vertices set $P'(\vec{v,w})
 [ \Sigma_{\vec{e} \ with \ i(\vec{e}) = v,
t(\vec{e}) = w} \ P(\vec{e}) ]$  be the sum over
all directed edges, say $\vec{e}$, with initial point $v$ and terminal
point $w$ of $P(\vec{e})$. On the other hand, if there is a directed
edge, say $\vec{\tau})$, with initial point $v$ and terminal vertex also $v$,
let  $\vec{loop}$ be a choice of  the direction of the associated loop
in $X'$  and set $P'(\vec{loop}) = (1/2) [ \Sigma_{\vec{e} \ with \ i(\vec{e}) = v,
t(\vec{e}) = v} \ P(\vec{e}) ]$. With these definitions, identifying
$C^0(X,P)$ with $C^0(X,P')$, since they have the same
vertices, one easily finds that the transmission graph Laplacians
are identical and so have the same eigenvalues.

Note however, the  ``Cheeger'' constants for $X,P$ and $X',P'$ appearing in
Theorem \ref{thmCheegerupper1} which estimate
the
same eigenvalue difference differ since they reflect differing geometry.

\vspace{.3in}
\section{A Categorical Approach. Push forward collapse.} \label{sectPushForward}

A more elegant  approach to the notion of transmission systems is as follows.

To each vertex $v$ of a graph $X$ assign a complex vector space, say $W(v)$.
To each directed edge, say $\vec{e}$ assign a linear mapping
$$
P(\vec{e}) : W( i(\vec{e})) \rightarrow W(f(\vec{e}))
$$

Now  define the $0$-cochains $C^0(X,P)$ as the set of mappings
$$
f: V(X) \rightarrow \bigsqcup_{v\in V(X) } \ W(v)
$$
where $\bigsqcup$ is the disjoint union so that
$f(v)$ lies in the vector space $W(v)$ associated to $v \in V(X)$.
[Often such $f$ are called sections.]

Now define a linear mapping
$$
A := C^0(X,P ) \rightarrow C^0(X,P)
$$
by
$$
(Af)(w) = \Sigma_{\vec{e} \ with \ i(\vec{e}) = w \ and \ t(\vec{e}) \ } \ P(\vec{e}) f(v)
$$

Let $T$ be the mapping $T^k : C^0(X,P ) \rightarrow C^0(X,P)$ defined by
$(Tf)(v) = (d_v)^k \ f(v)$.

Define the transmission graph Laplacian by
$$
\Delta^P = Id - T^{-1/2} A T^{-1/2} = T^{-1/2} ( T - A) T^{-1/2}
$$

As easily seen by taking bases, the over becomes identical
to the concrete representations given in terms of matrices,
where $A$ becomes the adjacency matrix $A^P$, etcetera.

\vspace{.3in}

Now if $\Psi : V(X) \rightarrow Y$ is a surjection of sets, then the
collapsed graph as explained in \S \ref{sectcollapse} has vertices
$Y$ and directed edges the same as those of $X$. An directed edge
$\vec{e}$ is to be thought of as having initial point
$\Psi( I(\vec{e}))$ and terminal point $\Psi(t(\vec{e}))$
in the collapsed graph. Call this collapsed graph $Y'$.

Now to each vertex, say $[p] \in Y$ associate the vector space
$$
W'([p]) = \bigoplus_{y \in V(x)  \ with \ \Psi(y) = p \ } \ W(y)
$$
the direct sum over the pre-images of $[p]$. This is usually called
the push forward of the vector spaces $\{ W(v)\}$.

Note that for $y \in V(x)$ there are the natural inclusions and
projections
$$
W(y) \stackrel{f_y}{\subset} \ W'(\Psi(y))  \stackrel{\pi_y}{\rightarrow} W(y)
$$
of the $y$-summand into the direct sum $W'(\Psi(y))$.

Define the push forward transmission system $P'$ on the collapsed
graph with vertices $Y$ by the ansatz: For directed edge $\vec{e}$ of $X$,
and so of $Y'$, define $P'(\vec{e})$ by the composite:

$$
\begin{array}{l}
P'(  \vec{e}  ) =        f_{ t(\vec{e})  }   \cdot  P(\vec{e}) \cdot \pi_{i(\vec{e})}  \\
W'(\Psi(i(\vec{e})))  = \bigoplus_{y \in V(x)  \ with \ \Psi(y) = \Psi(i(\vec{e})) \ } \ W(y) \\
     \rightarrow W'(\Psi(t(\vec{e}))) = \bigoplus_{z \in V(x)  \ with \ \Psi(z) = \Psi(t(\vec{e})) \ } \ W(z )\\
     \end{array}
     $$

The push forward of a   trivial transmission system,  that is one with all ranks equal to one
and all maps $P(\vec{e})$ equal the identity, become a system  with vector spaces of  varying ranks.
Possibly this is of interest in studying group actions on graphs.
It is to be noted the push forward transmission graph Laplacian
is just a reorganization of the earlier one, so has the same
eigenvalues and multiplicities. However, upon amalgamation
the resultant graphs are different so the ``Cheeger'' constants
of Theorem \ref{thmCheegerupper1}
are different.

\section{Morse Theory, Riemann Surfaces and Transmission systems.} \label{sectRiem}

An important case in which a graph arises naturally is in Morse theory.
If $M$ is a compact closed smooth Riemannian manifold, a Morse function
is a real valued function $f : M \rightarrow R$ with non-degenerate
critical points.
The critical points of the graph are the critical points of the Morse function
on the  manifold $M$
and the edges correspond to trajectories under the gradient flow
 which pass from the critical point of index $i$ to index $i-1$
 under the gradient flow. These arise as the intersection of
 the ascending and descending cells from this pair of critical points
 and are finite in number for a Morse function of Morse-Smale type.
 This directed and integer labeled graph when reinterpreted  gives
 a chain complex is a model which allows computation of the homology of $M$ \cite{Milnor1}

 In a similar way, if $E$ is a flat bundle over $M$, parallel  transport along
 these trajectories gives a invertible transmission system whose associated
 chain complex computes the homology of $M$ with coefficients in this local system.

 It is of great interest that Witten \cite{Witten1} has introduced a deformation of the
 Laplacian which has low lying eigen-modes localized at the critical points
 of  $f$ and whose tunneling is precisely along these descending trajectories. For the
 Laplacian coupled to a flat bundle as in \cite{CappellMiller}, this tunneling
 records the invertible transmission system obtained by parallel transport.

 \vspace{.3in}
 Another interesting example of a graph with natural transmission system arises
 in the context of Riemann surfaces.
Let a graph $X$ be embedded into a Riemann surface $S$
 so each edge and with its
end points is mapped smoothly.  The Riemannian metric on $S$ determines
 parallel transport along paths using the Riemannian connection.
 Now  to each vertex, say $v\in V(X)$, there is associated
the complex tangent space $T_vS$. Also one  may use parallel transport to
identify the tangent space $T_v$ at the initial point $v$ of an edge, say $\vec{e}$,
with the tangent space $T_w$ at the terminal point $w$  of the edge $\vec{e}$.
Since parallel transport preserves length, angle, and orientation of vectors,  this is a complex
linear mapping, which is denoted here by
$$
P(\vec{e}) : T_vS \rightarrow T_wS \ where \ v = i(\vec{e})\ and \ w = t(\vec{e})\ .
$$

In this way an invertible transmission system $P$ is naturally defined from the
embedding and the parallel transport over the Riemann surface. One might
imagine that the transmission graph Laplacian is in some crude
way the analogue of the $\overline{\partial}$-Laplacian of the Riemann surface
regarded as a complex Kahler variety.

\begin{flushleft}

Sylvain E. Cappell

Courant Institute, N.Y.U.

251 Mercer Street

New York, NY 10012

email: cappell@courant.nyu.edu

\vspace{.1in}

Edward Y. Miller

Mathematics Department

Polytechnic Institute of New York University

Six MetroTech Center

Brooklyn, NY 11201

email: emiller@poly.edu

\end{flushleft}

\end{document}